%% file: main.tex
\documentclass{ACmod}

\input{defs}

\usepackage{tikz}
\usepackage[tocflat]{tocstyle}
\usetikzlibrary{positioning}
\usetikzlibrary{cd}
\usepackage{amscd}
\usetikzlibrary{arrows,shapes,positioning}
\usetikzlibrary{decorations.markings}
\tikzstyle arrowstyle=[scale=1]
\tikzstyle directed=[postaction={decorate,decoration={markings,
    mark=at position .65 with {\arrow[arrowstyle]{stealth}}}}]
\tikzstyle reverse directed=[postaction={decorate,decoration={markings,
    mark=at position .65 with {\arrowreversed[arrowstyle]{stealth};}}}]
\usepackage[tocflat]{tocstyle}
\tikzset{cross/.style={cross out, draw=black, minimum size=2*(#1-\pgflinewidth), inner sep=0pt, outer sep=0pt},
cross/.default={1pt}}
\begin{document}

\input{title}

\input{abstract}

\maketitle

\setcounter{tocdepth}{1}
\tableofcontents

\input{intro}
\input{annuli}
\input{vertices}

\input{tcft}

\input{biblio}
\end{document}

%% file: defs.tex
\DeclareFontFamily{U}{rsf}{}
\DeclareFontShape{U}{rsf}{m}{n}{
  <5> <6> rsfs5 <7> <8> <9> rsfs7 <10-> rsfs10}{}
\DeclareMathAlphabet{\mathscr}{U}{rsf}{m}{n}

\DeclareMathAlphabet{\mathgth}{U}{euf}{m}{n}

\DeclareFontFamily{U}{cyr}{}
\DeclareFontShape{U}{cyr}{m}{n}{
  <5> wncyr5 <6> wncyr6 <7> wncyr7 <8> wncyr8 <9> wncyr9 <10-> wncyr10}{}
\DeclareMathAlphabet{\mathcyr}{U}{cyr}{m}{n}

\input cyracc.def

\DeclareSymbolFont{bbold}{U}{bbold}{m}{n}
\DeclareSymbolFontAlphabet{\mathbbold}{bbold}

\makeatletter
\def\operator@font{\sf}
\makeatother

\newcommand{\fr}{{\sf fr}}
\newcommand{\comb}{{\sf comb}}
\newcommand{\g}{\mathgth{g}}
\newcommand{\hg}{\widehat{\mathgth{g}}}
\newcommand{\h}{\mathgth{h}}
\newcommand{\hh}{\widehat{\mathgth{h}}}

\newcommand{\cA}{{\mathscr A}}

\newcommand{\cP}{{\mathscr P}}

\newcommand{\cR}{\mathcal{R}}

\newcommand{\cS}{{\mathscr S}}

\newcommand{\cV}{{\mathcal V}}

\newcommand{\bM}{{\overline{M}}}

\newcommand{\hcV}{\widehat{\cal V}}

\newcommand{\Annu}{{\mathcal{A}nn}^+}
\newcommand{\PEnd}{{\sf End}}
\newcommand{\cont}{{\sf c}}

\newcommand{\hS}{\mathsf{hS}}

\newcommand{\sgn}{{\mathsf{sgn}}}

\newcommand{\res}{{\mathsf{res}}}

\newcommand{\bbk}{\mathbb{K}}

\newcommand{\bbQ}{\mathbb{Q}}

\newcommand{\Tate}{\mathsf{Tate}}

\DeclareMathOperator{\Res}{Res}
\DeclareMathOperator{\Sym}{Sym}

\DeclareMathOperator{\id}{id}

\DeclareMathOperator{\Hom}{Hom}

\newcommand{\ra}{\rightarrow}

\newcommand{\lra}{\longrightarrow}

\newcommand{\Z}{\mathbb{Z}}

\newcommand{\chkX}{\check{X\,}\!}

\newcommand{\del}{\partial}

\renewcommand{\phi}{\varphi}

%% file: title.tex
\title{Categorical Enumerative Invariants, I: String vertices}

\author[C\u ald\u araru, Costello and Tu]{% 
Andrei C\u ald\u araru, Kevin Costello and Junwu Tu}

\address{Andrei C\u ald\u araru, Mathematics Department,
University of Wisconsin--Madison, 480 Lincoln Drive, Madison, WI
53706--1388, USA.}
\address{Kevin Costello, Perimeter Institute, Perimeter Institute of
Theoretical Physics, 31 Caroline St.\ N, Waterloo, ON N2L 2Y5, Canada}
\address{Junwu Tu, Institute of Mathematical Sciences, ShanghaiTech University, Shanghai, China, 201210.}

%% file: abstract.tex
\begin{abstract} 
  {\sc Abstract:} We define combinatorial counterparts
  to the geometric string vertices of Sen-Zwiebach and
  Costello-Zwiebach, which are certain closed subsets of the moduli
  spaces of curves.  Our combinatorial vertices contain the same
  information as the geometric ones, are effectively computable, and
  act on the Hochschild chains of a cyclic $A_\infty$-algebra.

  This is the first in a series of two papers where we define
  enumerative invariants associated to a pair consisting of a cyclic
  $A_\infty$-algebra and a splitting of the Hodge filtration on its
  cyclic homology.  These invariants conjecturally generalize the
  Gromov-Witten and Fan-Jarvis-Ruan-Witten invariants from symplectic
  geometry, and the Bershadsky-Cecotti-Ooguri-Vafa invariants
  from holomorphic geometry.
\end{abstract}

%% file: intro.tex
\section{Introduction}

\paragraph
This is the first in a series of two papers whose purpose is to define
enumerative invariants associated to a pair $(A,s)$ consisting of a
cyclic $A_\infty$-algebra $A$ and a splitting $s$ of the Hodge
filtration on its periodic cyclic homology.  

\noindent
These invariants conjecturally generalize
\begin{itemize}
\item the Gromov-Witten invariants of a symplectic space $X$ ($A$
  is obtained from the Fukaya category of $X$);
\item the Bershadsky-Cecotti-Ooguri-Vafa invariants of a complex
  Calabi-Yau manifold $\chkX$ ($A$ is obtained from the derived
  category of coherent sheaves on $\chkX$);
\item Saito's $g=0$ invariants associated to an isolated singularity
  ($A$ is obtained from the corresponding category of matrix
  factorizations);
\item the B-model analogue of the Fan-Jarvis-Ruan-Witten invariants
  (same as previous case).
\end{itemize}
In this paper we focus on defining combinatorial string vertices and
their action on Hochschild chains of $A$. These will be used in the
follow-up paper~\cite{CalTu} to give a constructive definition of the
categorical enumerative invariants, using the action of the Givental
group on the space of string field theories.

\paragraph
String vertices were introduced in physics by
Sen-Zwiebach~\cite{SenZwi} as an essential part of their construction
of string field theory.  The original approach defined string vertices
as certain subsets $\cV_{g,n}\subset M_{g,n}$ parametrizing curves
that satisfy a geometric condition in terms of minimal area metrics.
Another geometric construction of string vertices, this time using
hyperbolic geometry, was recently described by
Costello-Zwiebach~\cite{CosZwi}.

Intuitively, the string vertex $\cV_{g,n}\subset M_{g,n}$ is the
complement, in the compactified moduli space of curves, of an open
tubular neighborhood of the boundary.  For example, the string vertex
$\cV_{1,1}$ consists of the points in $M_{1,1}$ which are at least
$\varepsilon$ away (with respect to some metric on $\bM_{1,1}$) from the
point corresponding to the nodal curve.  As a singular chain
$\cV_{1,1}$ is not closed.  Its boundary (a real circle) can be
thought of as consisting of all elliptic curves that are obtained by
{\em twist-sewing} (see Section~\ref{sec:vertices}) on two boundary
circles of the unique surface in $M_{0,3}$.

\paragraph
This idea was expanded in~\cite{SenZwi} where it was noted that,
appropriately defined, the total string vertex
\[ \cV = \sum_{g,n\geq 1} \cV_{g,n} \, \hbar^g \lambda^{2g-2+n} \] 
should satisfy a form of the quantum master equation (QME)
\[ \del \cV + \hbar\Delta \cV + \frac{1}{2}\{\cV, \cV\} = 0 \]
in a Batalin-Vilkovisky (BV) algebra constructed from the spaces of
singular chains $C_*(M_{g,n})$.  Here, and in the sequel, $C_*$ will
denote the functor of normalized singular chains, and $\lambda$ and
$u$ will be formal variables of homological degree $-2$. The variable
$\hbar$ will have degree $-2$ in this setting, and degree $-2+2d$ in
the setting of an $A_\infty$-algebra (see Theorem B below).

This observation led Costello to a homological definition of string
vertices.  Let $\g$ denote the differential graded (dg) vector
space
\[ \g= \left (\bigoplus_{g,n} C_* (M_{g,n}^\fr)_{(S^1)^n \rtimes
      \Sigma_n}[1]\right )\series{\hbar, \lambda}, \]
where $M_{g,n}^\fr$ is the moduli space of curves of genus
$g$ with $n$ framed marked points.  It is an $(S^1)^n$ bundle over
the usual moduli space of curves $M_{g,n}$.

The space $\g$ carries geometric operations $\Delta$,
$\{\,-\,,\,-\,\}$ making it into (a shift of) the primitive part of a
BV algebra.  The geometrically-defined string vertices are a solution
to the QME.  Costello~\cite{Cos} shows that there
is a unique solution $\cV$ to the QME, up to homotopy,
once we fix the initial condition
\[ \cV_{0,3} = \frac{1}{6} \,[M_{0,3}^\fr]\,\lambda^1. \] 
Using this result Costello defines the string vertices $\cV_{g,n}$ as the
coefficients of the expansion of $\cV$ in terms of $g$ and $n$,
see~\cite[Theorem 1]{Cos}.   

Note that we do not need the full BV algebra structure
discussed above: all we need is the structure of differential graded
Lie algebra (DGLA) on $\g$ given by
\[ \g = \left (\g, \del + \hbar \Delta, \{\,-\,,\,-\,\} \right). \]
The string vertex $\cV$ is then the unique solution of the Maurer-Cartan
equation in this DGLA satisfying the initial condition above.

The goal of this paper and the sequel~\cite{CalTu} is to construct the
string vertices combinatorially, in the ribbon graph model for moduli
spaces. We will leverage the fact that string vertices are unique to
provide an algorithm to construct them explicitly.  The main goal of
this paper is to explain how to overcome an important technical
difficulty, which we now discuss.

\paragraph
Both the geometric construction of string vertices of Sen-Zwiebach and
the homological construction of Costello have two important drawbacks.
One is that it is not clear how to compute the string vertices
$\{\cV_{g,n}\}$ explicitly.  Another is that it is {\em a priori}
difficult to do computations in two-dimensional topological field
theory (2d TFT) using either version of string vertices.

The first problem comes from the fact that even though the dg vector
spaces $C_*(M_{g,n}^\fr)$ admit combinatorial descriptions in terms of
ribbon graphs, the operators $\Delta$, $\{\,-\,,\,-\,\}$ do not have
direct descriptions in these terms. (They do not respect the cell
decomposition of the moduli spaces of curves given by ribbon graphs.)
Hence the QME cannot be translated directly into an equation for
ribbon graphs.

The second issue is more conceptual.  In the process of constructing
enumerative invariants of a pair $(A,s)$ we need to let the string
vertices act on the space $L = CC_*(A)[d]$ of shifted Hochschild
chains of $A$. (Here $d$ denotes the Calabi-Yau dimension of $A$. 
With either the geometric or homological string vertices defined above
this is not possible.

The problem is that while the algebra $A$ gives rise to a {\em
  positive boundary} 2d TFT, described explicitly in terms of ribbon
graphs by Kontsevich-Soibelman~\cite{KonSoi} and Costello~\cite{Cos1},
the homological string vertices live in the zero-input part of the
TFT, so it is not clear how they act on $L$.  More explicitly, the
positive boundary 2d TFT obtained from $A$ is encoded in a PROP
action, given by degree zero operations
\[ \rho^A_{g,k,l}: C_*\left (M_{g,k,l}^\fr \right )[d\cdot(2-2g-2k)]
  \lra {{\sf Hom}}\left (L^{\otimes k}, L^{\otimes l} \right ) \] for
every $g,k\geq 1,l$.  Here $M_{g,k,l}^\fr$ denotes the moduli space of
Riemann surfaces of genus $g$ with $k$ incoming and $l$ outgoing
framed boundaries.  The homological string vertex $\cV_{g,n}$ is
naturally interpreted as a chain in $C_*(M_{g,0,n})$; the
positive boundary condition on the TFT defined from $A$ means that it
does not act on $L$ in any obvious way.

\paragraph
We note that Lurie~\cite{Lur} shows that one can remove the constraint
that $k \ge 1$ under the hypothesis that the algebra $A$ is smooth as
well as proper.  However, we do not know how to make Lurie's
construction \emph{explicit} in any combinatorial model for the moduli
of Riemann surfaces.

\paragraph
Both issues mentioned above arise from a common root problem, namely
that the calculations we want to carry out (finding string vertices,
and computing their action on Hochschild chains) are described in
terms of the spaces $C_*(M_{g,0,n})$.  These spaces do not have an
effective combinatorial description, even though for $k\geq 1$ the
spaces $C_*(M_{g,k,l})$ do.  The problems would be resolved if we had
a mechanism for replacing computations for $M_{g,0,n}$ by computations
for $M_{g,k,l}$ for $k\geq 1$.

The solution we propose is to use a resolution $\hg$ of the DGLA $\g$.
We shall call $\hg$ the {\em Koszul resolution} of $\g$.  Since
finding the string vertices amounts to solving the Maurer-Cartan
equation in $\g$, solving it in $\hg$ gives equivalent information.

\paragraph
First we need to introduce some notation.  Let $k$ and $l$ be non-negative
integers, let $n=k+l$, and let $V$ be a dg vector space endowed with
an action of $(S^1)^n \rtimes (\Sigma_k\times \Sigma_l)$. Examples of
such spaces include $C_*(M_{g,k,l}^\fr)$ and $L^{\otimes n}$.  (We
refer the reader to Section~\ref{sec:annuli} for details on circle
actions on dg vector spaces.)

The homotopy quotient of $V$ by the group
$(S^1)^n \rtimes (\Sigma_k \times \Sigma_l)$ will be denoted by
$V_\hS$.  This notation is ambiguous, but $k$ and $l$ will be evident
from context.  In this quotient we let $\Sigma_k$ and $\Sigma_l$ act
via the alternating and trivial representations, respectively.  If the
original differential on $V$ is $d$, we will denote the differential
on $V_\hS$ by $d+uB$.

\noindent
The following are the main theorems of this paper.

\paragraph {\bf Theorem A.}
{\em 
For any $g,k,l$ non-negative integers, including $k=0$, let $V_{g,k,l}
= C_*(M_{g,k,l}^\fr)_\hS$ and define
\begin{align*}
\g & = \left ( \bigoplus_{g, n\geq 1} V_{g,0,n}[1] \right ) \series{\hbar,
     \lambda}
\intertext{and}
\hg & = \left (\bigoplus_{g,k\geq 1, l} V_{g,k,l}[2-k]
      \right )\series{\hbar, \lambda}.
\end{align*}

Then there are operators
\begin{align*}
\Delta & : V_{g,k,l}\ra V_{g+1, k, l-2}[-1]\\
\{\,-\,,\,-\,\} & : V_{g,0,n} \otimes V_{g',0,n'} \ra V_{g+g', 0,n+n'-2}[-1] \\
\iota & : V_{g,k,l} \ra V_{g,k+1,l-1}\\
\{\,-\,,\,-\,\}_i & : V_{g,k,l}\otimes V_{g',k',l'} \ra
                    V_{g+g'+i-1,k+k'-i,l+l'-i}\,[-i]\mbox{ for
                    }k\geq 1, k'\geq 1, i\geq 1
\end{align*}
which endow $\g$ and $\hg$ with DGLA structures given by
\begin{align*}
  \g & = \left ( \g, \del + uB + \hbar\Delta, \{\,-\,,\,-\,\} \right )
  \intertext{and}
  \hg & = \left ( \hg, \del + uB+ \iota + \hbar \Delta, \{\,-\,,\,-\,\}_1 +
      \frac{1}{2!}\hbar\{\,-\,,\,-\,\}_2 + \frac{1}{3!}\hbar^2\{\,-\,,\,-\,\}_3 + \cdots
        \right ).
\end{align*}

The map $\iota:\g \ra \hg$ is a quasi-isomorphism of DGLAs.
Hence the Maurer-Cartan equation in $\hg$ admits a solution 
\[ \hcV = \sum_{g,k\geq 1,l} \hcV_{g,k,l}\, \hbar^g \lambda^{2g-2+k+l} \]
which is unique up to homotopy once we fix 
\[ \hcV_{0,1,2} = \iota\left (\cV_{0,3} \right )=
  \frac{1}{2}\,[M_{0,1,2}^\fr]\,\lambda^1. \] } 
\vspace{-1em}

\paragraph {\bf Remark.} The difference of a factor of three between the
initial conditions of the Maurer-Cartan equations in $\g$ and in $\hg$
is explained by the fact that $\iota$ relabels outputs into inputs,
and there are three outputs in $\cV_{0,3}$ to be relabeled.

\paragraph {\bf Theorem B.}
{\em 
Let $A$ be a smooth and proper cyclic $A_\infty$-algebra for which the
Hodge-de Rham degeneration property holds.  If $d$ is
the Calabi-Yau dimension of $A$ denote by $L = CC_*(A)[d]$ the
$d$-shifted Hochschild chain complex of $A$.  The variable $\hbar$
will be of degree $-2+2d$.

For any $k,l$ non-negative integers, including $k=0$, let
$W_{k,l}$ be the dg vector space 
\[ W_{k,l} = \Hom\left(L^{\otimes k}, L^{\otimes l}\right)_\hS. \]
Define dg vector spaces 
\begin{align*}
\h & = \left ( \bigoplus_{n}W_{0,n}[1-2d] \right ) \series{\hbar,
    \lambda} 
\intertext{and}
\hh & = \left (\bigoplus_{k\geq 1, l} W_{k,l}[2-2d+2dk-k] \right )\series{\hbar, \lambda}. 
\end{align*}

Then there are operators 
\begin{align*}
\Delta & :  W_{k,l} \ra W_{k, l-2}[2d-1]\\
\{\,-\,,\,-\,\} & :  W_{0,n} \otimes W_{0,n'}\ra W_{0,n+n'-2}[2d-1]\\
\iota & :  W_{k,l} \ra W_{k+1,l-1}[2d] \\
\{\,-\,,\,-\,\}_i & :  W_{k,l}\otimes W_{k',l'} \ra
                    W_{k+k'-i,l+l'-i}\,[-i] \mbox{ for } k\geq 1,
                    k'\geq 1, i\geq 1\\
\end{align*}
which endow $\h$ and $\hh$ with the structure of DGLAs
\begin{align*}
 \h & =  \left ( \h, b  + uB + \hbar \Delta, \{\,-\,,\,-\,\} \right ), \\
 \hh & =  \left ( \hh, b  + uB + \iota + \hbar \Delta, \{\,-\,,\,-\,\}_1 +
      \frac{1}{2!}\hbar\{\,-\,,\,-\,\}_2 + \frac{1}{3!}\hbar^2\{\,-\,,\,-\,\}_3 + \cdots
    \right ). 
\end{align*}
The map $\iota:\h \ra \hh$ is a quasi-isomorphism of DGLAs.  }

\paragraph
The DGLA structures on $\g$ and $\h$ are well known, having already
appeared in~\cite{Cos}.  These structures arise from the fact that the
spaces $V_{g,0,n}$ and $W_{0,n}$ admit actions of the {\em category}
of annuli: sewing the degree one annulus with one input and one output
gives the circle action, and sewing the one with two inputs induces
the $\Delta$ operator.

The new idea in this paper is that we have a bit more structure,
because the collection of annuli forms a {\em PROP} which acts.  Sewing
just one input of the degree zero annulus with two inputs defines the
$\iota$ operator.  It effectively relabels an output into an input, in
all possible ways; on Hochschild chains it acts by a formula that is
extremely similar to that of the Koszul differential.  This operator
is the crucial new ingredient in the definition of the resolutions
$\hg$ and $\hh$.

\paragraph
The DGLA $\hg$ and in particular the operators $\Delta$, $\iota$,
$\del$, $B$, $\{\,-\,,\,-\,\}_i$ admit combinatorial descriptions in
terms of ribbon graphs.  Thus the string vertex $\hcV$ can be computed
recursively as combinations of ribbon graphs, starting from the
initial condition.  Using the Kontsevich-Soibelman~\cite{KonSoi}
action $\rho^A$ (as described in detail in~\cite{CalChe}) we obtain
explicitly computable tensors
\[ \widehat{\alpha}^A_{g,k,l} = \rho^A(\hcV_{g,k,l}) \in W_{k,l} =
  \Hom(L^{\otimes k}, L^{\otimes l})_\hS = \Hom(\wedge^kL_+, \Sym^l
  L_-). \]
Here $L_+ = uL\series{u}$ and $L_- = L[u^{-1}]$.

As the image of the Maurer-Cartan element $\hcV$, the element
\[ \widehat{\alpha} = \sum_{g, k\geq 1, l} \widehat{\alpha}^A_{g,k,l}
  \hbar^g  \lambda^{2g-2+k+l} \]
also satisfies the Maurer-Cartan equation in $\hh$. (One can verify
carefully that it has degree $-1$ in $\hh$.)

\paragraph
\label{subsec:beta}
For the use in~\cite{CalTu} it will be more convenient to
replace the spaces
\[ W_{k,l} =  \Hom(\wedge^kL_+, \Sym^l L_-) \]
in the definition $\hh$ by the isomorphic ones (up to a shift by $-k$)
\[ W_{k,l}' = \Hom\left ( \Sym^k(L_+[1]), \Sym^l L_-\right ). \]
The resulting DGLA will still be denoted by $\hh$
\[ \hh = \left (\bigoplus_{k\geq 1, l} W_{k,l}'[2-2d+2dk] \right
  )\series{\hbar,\lambda}. \]
The old Maurer-Cartan element is identified with the element
\[ \widehat{\beta} = \sum_{g, k\geq 1, l} \widehat{\beta}^A_{g,k,l}
  \hbar^g  \lambda^{2g-2+k+l} \]
in the new setting.
The individual components $\widehat\beta^A_{g,k,l} \in \Hom\left (
  \Sym^k(L_+[1]), \Sym^l L_-\right )$ are related to the old ones
by the signs
\[\widehat{\beta}^A_{g,k,l}(x_1, \ldots, x_k) = (-1)^\epsilon\cdot
  \widehat{\alpha}^A_{g,k,l}(x_1, \ldots, x_k) \]
where
\[ \epsilon  = \sum_{i=1}^k (k-i)|x_i|. \]
The tensors $\widehat{\beta}^A_{g,k,l}$ will be one of the main
ingredients in the definition of the categorical invariants
in~\cite{CalTu}. 

\paragraph {\bf Standing assumptions.} We work over the field $\bbQ$
of rational numbers.  Whenever a pair $(g,n)$ or a triple $(g,k,l)$ of
non-negative integers appears (where $g$ refers to the genus of a
curve) it will be assumed that $2g-2+n>0$, $2g-2+k+l > 0$.  In any
summation over such triples, unconstrained variables are assumed to be
non-negative.  Thus $\sum_{g,n\geq 1}$ means the sum is over
all choices of $g\geq 0$, $n\geq 1$, $2g-2+n>0$.  

We use homological grading conventions: for a graded vector space
$V = \oplus V_n$, $V[k]$ is the graded vector space whose $n$-th
graded piece is $V_{n-k}$. Hence a degree zero map $V\ra W[k]$ {\em
  decreases} degree by $k$.

Outside of the introduction we will only work with the reduction of the
grading to $\Z/2\Z$.  In particular we will ignore even shifts (which
do not affect signs).  The degree shifts chosen in Theorems A and B have
been carefully designed to work for the $\Z$-grading.

In ribbon graph diagrams an input or output without a label is assumed
to have been labeled by $u^0$.  For a univalent white vertex we do
not mark the starting half-edge (since it is uniquely determined).

\paragraph {\bf Acknowledgments.} We would like to thank Nick
Sheridan, Si Li, and Dima Arinkin for patiently listening to the various
problems we ran into at different stages of the project, and for
providing insight.

Andrei C\u ald\u araru was partially supported by the National Science
Foundation through grant number DMS-1811925. The work of Kevin
Costello is supported by the NSERC Discovery Grant program and by the
Perimeter Institute for Theoretical Physics. Research at Perimeter
Institute is supported by the Government of Canada through Industry
Canada and by the Province of Ontario through the Ministry of Research
and Innovation. 

%% file: annuli.tex
\section{DGLAs induced by actions of the PROP of annuli}
\label{sec:annuli}

In this section we define the PROP of annuli $\Annu$ and we give two
different constructions of DGLAs associated to an action of $\Annu$
onto another PROP $\cP$.  The first construction yields the DGLAs
$\g$, $\h$ in the previous section; the second construction gives rise
to their resolutions $\hg$, $\hh$. 

Intuitively, both these constructions can be viewed as the results of
a two-part process.  First, the action of the generator $S$ of $\Annu$
gives rise to homological circle actions on the inputs and outputs of
the operations in $\cP(m,n)$.  The composition of $\cP$ does not
descend to the homotopy quotient $\cP_\hS$ by these circle actions,
but twisted sewing (composing with a circle action in between) does.
Second, the other generator $M$ of $\Annu$ induces two new structures
on the quotient $\cP_\hS$ -- twisted sewing of two outputs, and
relabeling of an output to an input.  These basic structures on
$\cP_\hS$ allow us to define the two types of DGLAs.

\paragraph
We will freely use the language of PROPs in this section. We refer
to~\cite{Mar} for generalities on PROPs. If $\cP$ is a PROP we use
$\otimes$ to denote horizontal composition:
\begin{align*}
\otimes & : \cP(n_1,m_1)\otimes\cdots\otimes \cP(n_s,m_s) \ra
           \cP(n_1+\cdots+n_s,m_1+\cdots+m_s) 
\intertext{and} 
\circ & : \cP(m,l)\otimes \cP(n,m) \ra \cP(n,l)
\end{align*}
the vertical composition. By definition, the space $\cP(n,m)$ is a
$\Sigma_n\times\Sigma_m$-bimodule. It is convenient to let $\Sigma_n$
act on the right, while $\Sigma_m$ will act on the left.

The unit element will be denoted by $\id\in \cP(1,1)$. For
$x\in \cP(n,m)$, $y\in \cP(m',l)$ we define the following notations:
\begin{itemize}
\item If $m'\leq m$, 
\[y\circ_{(i_1,\ldots,i_{m'})} x = (y\otimes
    \underbrace{\id\otimes\cdots\otimes\id}_{\text{$(m-m')$-copies}})\circ
    (\sigma \cdot x),\] where $(i_1,\ldots,i_{m'})$ is an ordered
  $m'$-tuple of distinct elements in $\{1,2,\cdots, m\}$, and
  $\sigma\in \Sigma_m$ is the permutation $(1,i_1)\cdots (m',i_{m'})$.
\item If $m'\geq m$, 
  \[ y_{(j_1,\ldots,j_m)} \circ x = (y\cdot \sigma) \circ
    (x\otimes\underbrace{\id\otimes\cdots\otimes\id}_{\text{$(m'-m)$-copies}}),\]
  where $(j_1,\ldots,j_{m})$ is an ordered $m$-tuple of distinct
  elements in $\{1,2,\cdots, m'\}$, and $\sigma\in \Sigma_{m'}$ is the
  permutation $(1,j_1)\cdots (m,j_{m})$. 
\end{itemize}

\paragraph{{\bf The PROP of annuli $\Annu$.}}
\label{para:annuli} 
The dg PROP $\Annu$ is defined by generators and relations as
follows. As a (unital) dg PROP it is generated by two operations:
\begin{enumerate}
\item A degree one element $S\in \Annu(1,1)$ such that $S\circ S=0$.
\item A degree zero element $M\in \Annu(2,0)$ which is invariant under
  the action of the symmetric group $\Sigma_2$.
\end{enumerate}
We put the zero differential on $\Annu$. The superscript ``$+$" is
used to
indicate that we only allow annuli with positive number of inputs.

Since $\Sigma_2$ acts trivially on $M$ we have 
\[M\circ_1 S= M\circ (S\otimes \id)=M\circ (\id\otimes S)=M\circ_2 S.\]
For this reason we shall simply write $\mathbb{M}=M\circ S$ for
either one of them. 

\paragraph
One can think of this PROP as the sub-PROP of the PROP of chains on
moduli spaces of curves which only includes annuli with positive
boundary.  Since this PROP is formal, we can replace it by its
homology, which is what we have described above.

\paragraph
According to the above description of $\Annu$, to give an
$\Annu$-algebra structure on a dg vector space $(V,b)$, i.e., a PROP
morphism 
\[ \lambda: \Annu \ra \PEnd(V),\]
is equivalent to giving the following data:

\begin{itemize}
\item A degree one map $B=\lambda(S): V\ra V$
  such that $B\circ B=0$, and $bB+Bb=0$.
\item A degree zero symmetric pairing
  $\langle\,-\,,\,-\,\rangle=\lambda(M): V\otimes V \ra \mathbb{C}$
  such that for any $x, y\in V$ we have
  \begin{align*}
    \langle b x, y\rangle +(-1)^{|x|}\langle x, b y\rangle &=0,\\
    \langle Bx, y\rangle - (-1)^{|x|}\langle x, By\rangle & = 0.
  \end{align*}
\end{itemize}
We emphasize that not all PROPs with $\Annu$ action are of this form,
only $\Annu$-algebras.  For example the PROP $\cS$ in the next section
is not of this type.

\paragraph{{\bf Circle actions.}}
We now recall some basic constructions regarding circle actions. Let
$X$ be a topological space endowed with an $(S^1)^n$-action. For
$1\leq i\leq n$ let $B_i: C_*(X) \rightarrow C_{*+1}(X)$ denote the
map that assigns to an $m$-chain $\sigma: \Delta^m \rightarrow X$ the
$m+1$-chain defined by
\[ S^1\times \Delta^m \stackrel{\id\times \sigma}{\longrightarrow}
  S^1\times X \rightarrow X,\]
with the second map given by the $i$-th circle action. Since we use
the normalized chain complex, we have $B_i^2=0$. We shall refer to
each operator $B_i$ as a homological circle action, or simply circle
action.

The $(S^1)^n$-equivariant homology of $X$ is computed by the chain
complex
\[ C_*(X)_{(S^1)^n} = \big( C_*(X)[u_1^{-1},\ldots,u_n^{-1}],
  \partial+\sum_{i=1}^n u_i \cdot B_i\big)\]
which we will denote by $C_*(X)_{(S^1)^n}$. There is a map of complexes
\[ \pi\circ \Res: C_*(X)_{(S^1)^n} \rightarrow
  C_*\big(X/(S^1)^n\big)\]
which takes the modified residue $\Res$ of an equivariant chain (the
coefficient of $(u_1\cdots u_n)^0$) and applies to it the canonical
projection map $\pi: C_*(X) \rightarrow C_*\left(
  X/(S^1)^n\right)$. If the $(S^1)^n$-action is free, then this map is
a quasi-isomorphism.

Note the slightly unusual convention where the equivariant homology is
computed by complexes starting with $u^0$, not $u^{-1}$.  Throughout
this paper the operator of taking residues means taking the
coefficient of $u^0$ and not of $u^{-1}$.  This is done in order to match
notation with~\cite{CalTu}.

\paragraph{\bf The twisted sewing operation.}
For the rest of this section we will place ourselves in the setting
where we have a dg PROP $\cP$ endowed with an $\Annu$ action; in other
words we have a morphism of PROPs $\lambda: \Annu\ra \cP$.

For any operation in $\cP(m,n)$ we will think of the operation of
composing with $\lambda(S)$, on either the inputs or the outputs, as
acting by a circle rotation.  (This is justified by the fact that
these compositions give rise to homological circle actions on
$\cP(m,n)$.)  Moreover, in addition to the usual composition of $\cP$
which sews an output of an operation to an input of another, the
$\Annu$ action also allows us to sew two outputs of the same operation
of $\cP$.  This is accomplished by composing these two outputs with
the two inputs of $\lambda(M)$.  We obtain a structure similar in
nature to that of a modular operad (though not exactly the same,
because we do not allow sewing of two inputs).

Our goal is to pass to the quotient $\cP_\hS$ of $\cP$ by the circle actions
induced from $\Annu$, in such a way that the two sewing operations
described above would descend to the quotient.  We take our cue from
the theory of moduli spaces of curves.  Sewing operations are
well defined on $M_{g,n}^\fr$, but after passing to the quotient
$M_{g,n}$ of $M_{g,n}^\fr$ by the circle actions these operations are
no longer well defined.  However, there is a modification of these
sewing operations that does descend, at least at the level of
singular chains: the {\em twisted sewing} operation, which first
performs an $S^1$ twist before sewing. 

This suggests that on the quotient $\cP_\hS$ there will be two fundamental
types of new compositions which we will call twisted sewings.  The
twisted self-sewing $\Delta$ is obtained by sewing two outputs of
an operation with $\lambda(\mathbb{M})$.  The twisted sewing between two
operations is obtained by composing an output of one with
a circle action $\lambda(S)$ and then with an input of the other.

These twisted sewings are the building blocks of the DGLAs that we
construct in the rest of this section from a PROP with $\Annu$ action.

\paragraph{{\bf A first DGLA associated to a PROP with $\Annu$ action.}}
We now make the above ideas precise.  Let $\cP$ be a PROP with $\Annu$
action, and consider the space $\cP(0,n)$ of operations with zero
inputs and $n$ outputs. For each $1\leq j\leq n$ we obtain a
homological circle action on $\cP(0,n)$ by composing with $\lambda(S)$
at the $j$-th output. Denote this circle action by
\[ B_j: \cP(0,n) \ra \cP(0,n), \;\;\; B_j(x)=\lambda(S)\circ_j x.\]
The associated $(S^1)^n$-equivariant chain complex is of the form
\[ \cP(0,n)_{(S^1)^n} = \big(\cP(0,n)
  [u_1^{-1},\ldots,u_n^{-1}], \partial+\sum_{i=1}^n u_i \cdot
  B_i\big),\]
where $\partial$ is the boundary map of $\cP(0,n)$. The symmetric group
$\Sigma_n$ still acts by the PROP structure on $\cP(0,n)_{(S^1)^n}$ while
also permuting the indices of the circle parameters $u_1,\ldots,
u_n$. Following our conventions we denote by
$\cP(0,n)_\hS$ the homotopy quotient of $\cP(0,n)$ by the
semidirect product $(S^1)^n\rtimes \Sigma_n$-action.

We shall define a DGLA structure on the graded vector space
\[ \mathfrak{g}^{\cP}= \left (\bigoplus_{n\geq 0} \cP(0,n)_\hS[1]\right )\series{\hbar}.\]
First, for an element
\[ x=\sum_{k_1,\ldots,k_n\geq 1} x_{k_1\ldots k_n} u_1^{-k_1}\cdots
u_n^{-k_n} \in \cP(0,n)
[u_1^{-1},\ldots,u_n^{-1}] \]
we set
\begin{align*}
\Delta (x) & =\sum_{1\leq i<j\leq n} \Res_{u_i=0,u_j=0} \lambda(\mathbb{M})\circ_{(i,j)} x\\
&= \sum_{\substack{1\leq i<j\leq n\\k_i=k_j = 0}} \lambda(\mathbb{M})\circ_{(i,j)}x_{k_1\ldots k_n} u_1^{-k_1}\cdots \widehat{u_i^0}\cdots \widehat{u_j^0}\cdots u_n^{-k_n}.
\end{align*}
In other words the operator $\Delta$ performs a twisted sewing
operation on all pairs of outputs of $x$ that are labeled by
$u^0$. The operator $\Delta$ is $\Sigma_n$-invariant, hence
it induces a map  on quotients which we still denote by $\Delta$,
\[ \Delta: \cP(0,n)_\hS\ra \cP(0,n-2)_\hS[-1].\]

\begin{Lemma}
\label{lem:differential}
The sum $\partial+\sum_{i=1}^n u_i \cdot B_i+\hbar\Delta$ is a
differential, i.e. we have
$(\partial+\sum_{i=1}^n u_i \cdot B_i+\hbar\Delta)^2=0$.
\end{Lemma}
\bigskip

\begin{Proof}
  We need to prove that $[\partial,\Delta]=0$, $[B_i,\Delta]=0$, and
  $\Delta^2=0$. The first identity follows from the fact that
  $\mathbb{M}$ is closed. Since
  $(M\circ S)\circ S=M\circ (S\circ S)=0$ the second identity
  holds. The last identity follows from the fact that the degree of
  $\mathbb{M}$ is odd.
\end{Proof}

\paragraph
Next we define a Lie bracket on $\mathfrak{g}^{\cP}$. The definition
is very similar to that of $\Delta$ above: we perform the twisted
sewing operation between an output of one element and an output of
another one.  More formally, for two elements of the form
\begin{align*}
x & = \alpha \cdot u_1^{-k_1}\cdots u_n^{-k_n} \in \cP(0,n) [u_1^{-1},\ldots,u_n^{-1}], \\
y & = \beta\cdot u_1^{-l_1}\cdots u_m^{-l_m} \in \cP(0,m) [u_1^{-1},\ldots,u_m^{-1}]\\
\end{align*}
we define their Lie bracket by
\[ \{x,y\} = (-1)^{|x|}\sum_{k_i=l_j=0}
  \lambda(\mathbb{M})\circ_{(i,n+j)} (\alpha\otimes\beta)
  u_1^{-k_1}\cdots\widehat{u_i^0}\cdots u_n^{-k_n}
  u_{n+1}^{-l_1}\cdots\widehat{u_{n+j}^0}\cdots u_{n+m}^{-l_m}.\]
The sum is over all the pairs $(i,j)$ such that $k_i=l_j=0$.  We
extend the bracket to all elements in $\g^\cP$ by linearity.  The
result is $\Sigma_n\times\Sigma_m$-invariant, hence it induces a Lie
bracket of degree one
\[ \{ -,- \}: \cP(0,n)_\hS\otimes \cP(0,m)_\hS \ra
  \cP(0,n+m-2)_\hS[-1].\]

\begin{Theorem}
\label{thm:dgla1}
The triple $\left( \mathfrak{g}^{\cP}, \partial+uB+\hbar\Delta,
  \{-,-\}\right)$ forms a DGLA. 
\end{Theorem}

\begin{Proof}
  We shall only prove the Jacobi identity. The Leibniz rule can be
  proved similarly. To avoid tedious formulas it will be useful to
  depict $\{x,y\}$ as
\[\{x,y\}=\begin{tikzpicture}[box/.style={draw,rounded corners,text width=.5cm,align=center},scale=0.3,baseline={(current bounding box.center)}]
\node[box] at (-2.5,0) (b) {$x$};
\node[box] at (2.5,0) (c) {$y$};
\draw [thick] (b) to [out=-90, in=-90] (c);
\end{tikzpicture}\]
where we put $\lambda(\mathbb{M})$ on the connecting arc between $x$
and $y$. With this notation the terms in the Jacobi identity can be depicted as
\begin{align*}
\{x,\{y,z\}\} & = \begin{tikzpicture}[box/.style={draw,rounded corners,text width=.5cm,align=center},scale=0.3,baseline={(current bounding box.center)}]
\node[box] at (-2,0) (b) {$x$};
\node[box] at (2,0) (c) {$y$};
\node[box] at (6,0) (d) {$z$};
\draw [thick] (b) to [out=-90, in=-90] node[midway,below]{} (1.8,-1);
\draw [thick] (2.2,-1) to [out=-90, in=-90] node[midway,below]{} (d);
\end{tikzpicture}+ \begin{tikzpicture}[box/.style={draw,rounded corners,text width=.5cm,align=center},scale=0.3,baseline={(current bounding box.center)}]
\node[box] at (-2,0) (b) {$x$};
\node[box] at (2,0) (c) {$y$};
\node[box] at (6,0) (d) {$z$};
\draw [thick] (b) to [out=-90, in=-90] node[midway,below]{} (6.4,-.8);
\draw [thick] (2.2,-1) to [out=-90, in=-90] node[midway,below]{} (d);
\end{tikzpicture}\\
\{\{x,y\},z\} & = \begin{tikzpicture}[box/.style={draw,rounded corners,text width=.5cm,align=center},scale=0.3,baseline={(current bounding box.center)}]
\node[box] at (-2,0) (b) {$x$};
\node[box] at (2,0) (c) {$y$};
\node[box] at (6,0) (d) {$z$};
\draw [thick] (b) to [out=-90, in=-90] node[midway,below]{} (1.8,-1);
\draw [thick] (2.2,-1) to [out=-90, in=-90] node[midway,below]{} (d);
\end{tikzpicture}+ \begin{tikzpicture}[box/.style={draw,rounded corners,text width=.5cm,align=center},scale=0.3,baseline={(current bounding box.center)}]
\node[box] at (-2,0) (b) {$x$};
\node[box] at (2,0) (c) {$y$};
\node[box] at (6,0) (d) {$z$};
\draw [thick] (-1.5,-.8) to [out=-90, in=-90] node[midway,below]{} (1.8,-1);
\draw [thick] (-2.5,-.8) to [out=-90, in=-90] node[midway,below]{} (6.4,-.8);
\end{tikzpicture}\\
\{y,\{x,z\}\}& = \begin{tikzpicture}[box/.style={draw,rounded corners,text width=.5cm,align=center},scale=0.3,baseline={(current bounding box.center)}]
\node[box] at (-2,0) (b) {$y$};
\node[box] at (2,0) (c) {$x$};
\node[box] at (6,0) (d) {$z$};
\draw [thick] (b) to [out=-90, in=-90] node[midway,below]{} (1.8,-1);
\draw [thick] (2.2,-1) to [out=-90, in=-90] node[midway,below]{} (d);
\end{tikzpicture}+ \begin{tikzpicture}[box/.style={draw,rounded corners,text width=.5cm,align=center},scale=0.3,baseline={(current bounding box.center)}]
\node[box] at (-2,0) (b) {$y$};
\node[box] at (2,0) (c) {$x$};
\node[box] at (6,0) (d) {$z$};
\draw [thick] (b) to [out=-90, in=-90] node[midway,below]{} (6.4,-.8);
\draw [thick] (2.2,-1) to [out=-90, in=-90] node[midway,below]{} (d);
\end{tikzpicture}\\
\end{align*}
\vspace*{-13mm}

Let us illustrate with one of the terms above how the signs work
out. The first term in $\{x,\{y,z\}\}$ is given by 
\begin{align*}
 & (-1)^{|x|+|y|} \lambda(\mathbb{M})\circ_{(i,j)} \big( x\otimes \lambda(\mathbb{M})\circ_{(k,l)} (y\otimes z)\big)\\
 = & (-1)^{|y|+1}  \lambda(\mathbb{M})\circ_{(k,l)} \big( \lambda(\mathbb{M})\circ_{(i,j)}(x\otimes y)\otimes z\big)\\
=  & (-1)^{|x|+|y|+1+|x|}  \lambda(\mathbb{M})\circ_{(k,l)} \big( \lambda(\mathbb{M})\circ_{(i,j)}(x\otimes y)\otimes z\big)
 \end{align*}
 The last expression is precisely the first term in
 $\{\{x,y\},z\}$. The other two signs can be checked similarly.
\end{Proof}

\paragraph{{\bf A DGLA associated with an $\Annu$-algebra.}}
\label{para:example1}
Let $V$ be an $\Annu$-algebra given by a dg PROP morphism
$\lambda: \Annu\ra \PEnd(V)$. As an example of the previous
construction we describe explicitly the DGLA $\mathfrak{g}^V$. We
denote by $b$ the internal differential of the dg vector space $V$, by
$B$ the image of $S$ under $\lambda$, and by
$\langle\,-\,,\,-\,\rangle$ the image of $M$.

Using the circle action $B$ form the Tate complex
\[ V^{\sf Tate}=\big( V\laurent{u}, b+uB \big).\] 
Denote by $V_+$ the subcomplex $\big( uV\series{u}, b+uB \big)$ of
$V^{\sf Tate}$, and by $V_-$ the quotient complex
$V^{\sf Tate}/V_+=\big( V[u^{-1}], b+uB\big)$.  

As a graded vector space the DGLA $\mathfrak{g}^V$ constructed by
Theorem~\ref{thm:dgla1} is given by
\[ \mathfrak{g}^V= \left ( (\Sym V_-) [1]\right )\series{\hbar}.\]
Define a pairing
 \[ \Omega:V_-\otimes V_- \rightarrow \mathbb{C}\]
which maps $x=u^{0}x_{0}+\cdots$ and $y=u^{0} y_{0}+\cdots$ to 
\[ \Omega(x,y)=\langle Bx_{0},y_{0}\rangle.\] 
The pairing $\Omega$ is antisymmetric because $B$ is self-adjoint
with respect to the pairing, which is itself symmetric. One can also
check that $\Omega$ is a chain map.  Since $B$ has degree $1$ and
$\langle-,-\rangle$ has degree $0$, the degree of $\Omega$ is $1$.

The twisted self-sewing operator $\Delta$ acts on
$ \mathfrak{g}^V= (\Sym V_- )[1]\series{\hbar}$ by
\begin{equation*}
 \Delta(x_1\cdots x_n)=\sum_{1\leq i<j\leq n} (-1)^\epsilon \Omega(x_i,x_j)x_1\cdots\widehat{x_i}\cdots\widehat{x_j}\cdots x_n.
 \end{equation*}
The sign 
\[ (-1)^\epsilon =
  (-1)^{|x_1|+\cdots+|x_{i-1}|+(|x_i|+1)(|x_{i+1}|+\cdots+|x_{j-1}|)} \]
is simply the Koszul sign obtained from passing $\Omega$ (an odd
operator) past the first $i-1$ elements, and then the pair $(\Omega,
x_i)$ past the next $j-i-1$ elements in the sequence $x_1 x_2\cdots x_n$.

The operator $\Delta$ is a BV operator (i.e., a second order
differential operator of homological degree and square zero) with
respect to the product structure on $\Sym V_-$.  The difference
$\Delta(x\cdot y) -(\Delta x \cdot y)-(-1)^{|x|}(x\cdot \Delta y)$,
which measures the failure of $\Delta$ to be a derivation, is a
symmetric operation in $x$ and $y$ of degree $1$. On the shifted
symmetric product $(\Sym V_-)[1]$ this failure is then an
anti-symmetric operation of degree zero which is precisely the Lie
bracket of the DGLA $\g^V$
\[ \{x,y\}= (-1)^{|x|}\big( \Delta(x\cdot y) -(\Delta x \cdot
  y)-(-1)^{|x|}(x\cdot \Delta y)\big).\]

\paragraph
The DGLA $\g^\cP$ associated to an $\Annu$-algebra $\cP$ can be
concretely described in an entirely similar way.  For each $n$ let
$\cP(0,n)^\Tate$ be the Tate complex for the action of $(S^1)^n$: if
$b$ is the differential on $\cP(0,n)$ and $B_i$ is the
$i$-th circle action, then
\[ \cP(0,n)^\Tate = \left ( \cP(0,n)\laurent{u_1,..,u_n}, b + \sum u_i
    B_i \right ). \]
We define $\cP(0,n)_-$ to be the quotient of $\cP(0,n)^\Tate$ by the
subcomplex spanned topologically by expressions of the form
$\alpha \cdot u_1^{k_1} \dots u_n^{k_n}$ with $\alpha\in\cP(0,n)$ and at
least one $k_i$ is positive.

The operator $\Omega$ defined above has an analog in this situation,
which is a map $\Omega_{ij} : \cP(0,n)_- \to \cP(0,n-2)_-$ for each
$1 \le i \neq j \le n$.

As before, we define $\cP(0,n)_\hS$ to be the coinvariants of
$\cP(0,n)_-$ by the symmetric group action (which acts simultaneously
on the $u_i$ variables).  We build a BV algebra structure on the
direct sum of the spaces $\cP(0,n)_\hS$.  The BV operator
$\Delta$ is the symmetrization of the operators $\Omega_{ij}$ and the
product comes from the horizontal product in the PROP structure.  With
the induced bracket $\{\,-\,,\,-\,\}$ (which measures the failure
of $\Delta$ to be a derivation) this yields a DGLA 
\[ \g^\cP = \left ( \bigoplus \cP(0,n)_\hS, b+\hbar \Delta, \{\,-\,,\,-\,\} \right ). \]

\paragraph{{\bf A second construction of a DGLA from an $\Annu$
    action.}} 
The previous construction built a DGLA from the zero-input part of a
PROP $\cP$ with $\Annu$ action.  We will now construct a different
DGLA from the positive boundary part of the same PROP.

A dg PROP $\cP$ is said to have positive boundary if $\cP(0,n)=0$ for
all $n\geq 1$. If $\cP$ is any dg PROP, its truncation $\cP^+$ to the
part with positive inputs is a positive boundary PROP. The PROP of
annuli $\Annu$ is also a positive PROP.

Let $\cP$ be a positive boundary PROP, and let $\lambda: \Annu\ra \cP$
be a morphism of dg PROPs. As before, composing with $\lambda(S)$
along the inputs or outputs defines $k+l$ commuting homological circle
actions on the space $\cP(k,l)$. The $(S^1)^{k+l}$-equivariant
chain complex is
\[ \cP(k,l)_{(S^1)^{k+1}} = \big( \cP(k,l)
  [u_1^{-1},\ldots,u_{k+l}^{-1}], \partial+\sum_{i=1}^{k+l} u_i \cdot
  B_i\big).\] 
Both $\Sigma_k$ and $\Sigma_l$ act on this chain complex.  We modify
the action of the symmetric group $\Sigma_k$ on the inputs of the
operations in $\cP(k,l)$ by twisting it by the sign representation
${\sf sgn}_k$.  We leave the action of $\Sigma_l$ on the outputs
unchanged.  As before, these groups act on the variables
$u_1,\ldots, u_k$, $u_{k+1},\ldots, u_{k+l}$ by permuting them.

There are two reasons behind the presence of the ${{\sf sgn}}_k$
representation in the action on inputs.  One is the presence of the
shift by $2-k$ of $\cP(k,l)$: it can be thought of as having forced
the $k$ inputs to be odd.  Another is related to the fact that we want
to obtain an analogue of the Koszul complex; the inputs will play the
role of exterior tensors, while the outputs will play the role of
symmetric tensors. Again, we shall follow our convention to use the
notation $\cP(k,l)_\hS$ to denote the homotopy quotient of
$\cP(k,l)$ by the action of the semidirect product
$(S^1)^{k+l}\rtimes (\Sigma_k \times \Sigma_l)$.

We shall construct a DGLA structure on the dg vector space
\[ \hg^\cP = \left ( \bigoplus_{k\geq 1, l} \cP(k,l)_\hS[2-k]\right )\series{\hbar}.\] 
Due to the shifts of degree in the definition of
$\widehat{\mathfrak{g}}^\cP$ it is necessary to deal with signs before
moving forward.  Indeed, the degree of an element
$x\in \cP(k,l)_\hS[2-k]$ is given by
\[ |x|=\deg(x)+2-k,\] 
where $\deg(x)$ stands for the chain degree of
$x\in \cP(k,l)_\hS$. On the shifted complex
$\cP(k,l)_\hS[2-k]$ the differential is given by
\[ \partial'+\sum_{i=1}^{k+l} u_i \cdot B_i' = (-1)^k\cdot \big( \partial+\sum_{i=1}^{k+l} u_i \cdot B_i\big).\]
where the notation $'$ is used to denote an operator after performing
the appropriate shifts of degrees. It is also most convenient to work
with a shifted version of the original PROP composition. For two elements
$x\in \cP(m,l)$, $y\in \cP(n,m)$, and two ordered indices
$I=(i_1,\cdots,i_r)$ and $J=(j_1,\cdots,j_r)$ we define the shifted
composition to be
\[ x_{\; I}\circ'_J y = (-1)^{(n+r)\cdot \deg(x)} x_{\; I}\circ_J y  \in \cP(n+m-r,l+m-r).\]

\noindent
The following lemma can be verified directly from the definitions.

\begin{Lemma}
\label{lem:leibniz}
With notations as above the following identities hold:
\begin{enumerate}
\item $\partial' (x_{\; I}\circ'_J y) = (\partial' x)_{\; I}\circ'_J y + (-1)^{r+|x|} x_{\; I}\circ'_J (\partial' y)$.
\item $(x_{\; I}\circ'_J y)_{\; K} \circ'_L z= x_{\; I}\circ'_J (y_{\; K}\circ'_L z)$.
\end{enumerate}
\end{Lemma}
\bigskip

\paragraph
We may now proceed to the construction of a DGLA structure on
$\widehat{\mathfrak{g}}^\cP$. First we shall construct a deformed
differential as follows. In addition to the existing
differential $\partial'+\sum_{i=1}^{k+l} u_i \cdot B'_i$, there are
two more components. One of them is similar to the previously defined
twisted self-sewing operation $\Delta$: we set
\[\Delta (x) =\sum_{1\leq i<j\leq l} \Res_{u_i=0,u_j=0}
  \lambda(\mathbb{M})\circ'_{(i,j)} x.\]
In other words we perform the twisted self-sewing operation only
among any pair of outputs that are labeled by $0$-powers of the
$u$'s. 

The other component of the differential is similar to the usual Koszul
differential: it changes an output of an operation $x$ to an input by
sewing that output with one of the inputs of $\lambda(M)\in \cP(2,0)$.
Formally we define
\begin{align*}
 \iota&: \cP(k,l)_\hS \ra \cP(k+1,l-1)_\hS\\
 \iota(x) &=  \sum_{1\leq i\leq l} \lambda(M)\circ'_i x.
 \end{align*}

\begin{Lemma}
The following identity holds
\[ \big( \partial'+\sum_{i=1}^{k+l} u_i \cdot B'_i+\iota+\hbar\Delta\big)^2=0.\]
\end{Lemma}

\begin{Proof}
  This follows from Lemma~\ref{lem:leibniz}. For example, to prove
  that $\iota\partial'+\partial'\iota=0$, we have

\begin{align*}
 \partial'\iota(x)  & = \sum_i \partial' (\lambda(M)\circ'_i x) \\
& =  \sum_i - \lambda(M)\circ'_i \partial' x  \;\; (\mbox{by Lemma~\ref{lem:leibniz}})\\
&  =  -\iota \partial' x.
 \end{align*}
 Note that when applying Lemma~\ref{lem:leibniz} we have used the facts that
 $|\lambda(M)|=0$ and $r=1$.  The other commutator identities can be
 verified in a similar way.
\end{Proof}

\paragraph
The construction of a Lie bracket on $\widehat{\mathfrak{g}}^\cP$ is
more involved, due to the asymmetry between inputs and outputs. Here
are two simple observations:
\begin{itemize}
\item[--] To keep the number of inputs positive we should never perform
  the twisted sewing of two inputs, otherwise the Lie bracket between
  two elements in $\cP(1,l_1)$ and $\cP(1,l_2)$ would produce elements
  with zero inputs.
\item[--] Since we want the bracket on $\hg^\cP$ to have degree zero, we
  should not perform twisted sewing of two outputs (which would be a
  degree $-1$ operation).
\end{itemize}

\paragraph
We conclude that in the definition of the bracket we should only ever
perform twisted sewings between inputs and outputs. As a first
attempt, for elements $x\in\cP(l,k)_\hS$ and
$y\in \cP(n,m)_\hS$, set
\[ \{x, y\}_1= \sum_{\substack{1\leq i \leq l\\1\leq j\leq m}}
  \Res (x {\,}_{i}\circ' \lambda(S) \circ'_{j} y) -
  (-1)^{|x||y|}\sum_{\substack{1\leq i \leq k\\1\leq j\leq n}}
  \Res (y {\,}_{j}\circ' \lambda(S) \circ'_{i} x)\] where
the residue is taken at $u_{i}=u_{j}=0$. 

It turns out that this operation is indeed a Lie bracket; however,
$\Delta$ fails to be a derivation of the bracket $\{-,-\}_1$, so we
don't yet get a DGLA:
\[ \Delta\{x,y\}_1-\{\Delta x, y\}_1- (-1)^{|x|}\{x,\Delta y\}_1\neq
  0.\] 
In terms of diagrams as in the proof of Theorem~\ref{thm:dgla1} this
discrepancy can be understood as the terms below:
\[\begin{tikzpicture}[box/.style={draw,rounded corners,text width=.5cm,align=center},scale=0.3,baseline={(current bounding box.center)}]
\node[box] at (-2.5,0) (b) {$x$};
\node[box] at (2.5,0) (c) {$y$};
\draw [thick,->] (b) to [out=-90, in=-90] (c);
\draw [thick] (-3.2,-.8) to [out=-90, in=-90] (3.2,-1);
\end{tikzpicture}+\begin{tikzpicture}[box/.style={draw,rounded corners,text width=.5cm,align=center},scale=0.3,baseline={(current bounding box.center)}]
\node[box] at (-2.5,0) (b) {$x$};
\node[box] at (2.5,0) (c) {$y$};
\draw [thick,<-] (b) to [out=-90, in=-90] (c);
\draw [thick] (-3.2,-.8) to [out=-90, in=-90] (3.2,-1);
\end{tikzpicture}.\]
Here the connecting arc with no arrow is the result of applying the
operator $\Delta$, while the directed arc is what is obtained by
applying the Lie bracket $\{-,-\}_1$.

To remedy this problem we introduce a new bracket operator $\{-,-\}_2$,
represented pictorially by
\[ \{x,y\}_2= \begin{tikzpicture}[box/.style={draw,rounded corners,text width=.5cm,align=center},scale=0.3,baseline={(current bounding box.center)}]
\node[box] at (-2.5,0) (b) {$x$};
\node[box] at (2.5,0) (c) {$y$};
\draw [thick,->] (b) to [out=-90, in=-90] (c);
\draw [thick,->] (-3.2,-.8) to [out=-90, in=-90] (3.2,-1);
\end{tikzpicture}+\begin{tikzpicture}[box/.style={draw,rounded corners,text width=.5cm,align=center},scale=0.3,baseline={(current bounding box.center)}]
\node[box] at (-2.5,0) (b) {$x$};
\node[box] at (2.5,0) (c) {$y$};
\draw [thick,<-] (b) to [out=-90, in=-90] (c);
\draw [thick,<-] (-3.2,-.8) to [out=-90, in=-90] (3.2,-1);
\end{tikzpicture}.\] 
The operator $\iota$ fails to be a derivation of $\{-,-\}_2$, and its
failure cancels precisely
$\Delta\{x,y\}_1-\{\Delta x, y\}_1- (-1)^{|x|}\{x,\Delta y\}_1$. But
then $\{-,-\}_2$ again is not compatible with the operator
$\Delta$. Thus we need to introduce a third bracket $\{-,-\}_3$, and
so on.

\paragraph
More systematically, for each $r\geq 1$ we introduce a bracket
operation which performs twisted sewing of any $r$ inputs with any $r$
outputs.  Specifically, for $x\in \cP(l,k)_\hS$ and $y\in \cP(n,m)_\hS$ we define
\[ \{x,
  y\}_r=\sum_{\substack{I=(i_1,\ldots,i_r),\\J=(j_1,\ldots,j_r)}}
  \Res(x {\,}_I\circ' \lambda(S)^{\otimes r} \circ'_{J} y)
  -(-1)^{|x||y|}
  \sum_{\substack{I=(i_1,\ldots,i_r),\\J=(j_1,\ldots,j_r)}} \Res (y {\,}_J\circ' \lambda(S)^{\otimes r} \circ'_{I} x)\]
where the residue is taken at $u_{I}=0, u_{J}=0$, and in the first sum 
$I=(i_1,\ldots,i_r)$ and $J=(j_1,\ldots,j_r)$ are $r$-tuples of
distinct elements in the sets
$\{1,\ldots,l\}$, $\{1,\ldots,m\}$, respectively.
The $r$-tuples $I$ and $J$ in the second summation are defined
similarly, but they are subsets of $\{1,\ldots,k\}$, $\{1,\ldots,n\}$.

\begin{Theorem}
\label{thm:dgla2}
The following graded vector space forms a DGLA
\[ \g^\cP = \left (\bigoplus_{k\geq 1, l} \cP(k,l)_\hS[2-k] \right )
  \series{\hbar} \] 
when endowed with
\begin{itemize}
\item[--] differential given by $\partial'+\sum_{i=1}^{k+l} u_i \cdot B'_i+\iota+\hbar\Delta$,
\item[--] Lie bracket given by $\{-,-\}_\hbar=\displaystyle\sum_{r\geq 1} \frac{1}{r!}\{-,-\}_r\cdot\hbar^{r-1}$.
\end{itemize}
\end{Theorem}

\begin{Proof}
For the Leibniz rule we need to prove that
\begin{align*}
(r+1)\big( \Delta\{x,y\}_r- & \{\Delta x, y\}_r-(-1)^{|x|}\{x,\Delta
  y\}_r \big) = \\
& = -\iota\{x,y\}_{r+1}+\{\iota x,y\}_{r+1} +(-1)^{|x|}\{x, \iota y\}_{r+1}.
\end{align*}
Let us begin with the right hand side. The terms in $-\iota\{x,y\}_{r+1}$
are of the form
\[ -\lambda(M)\circ'_p \Res(x_{\; I'}\circ'
  \lambda(S)^{r+1}\circ'_{J'} y)+(-1)^{|x||y|} \lambda(M) \circ'_q
  \Res (y_{\; K'} \circ' \lambda(S)^{r+1} \circ'_{L'} x),\] 
with the multi-indices $I'$, $J'$, $K'$, and $L'$ of cardinality
$r+1$. This further implies to the following
\begin{align*}
\begin{split}
&\underbrace{ -\Res\big( (\lambda(M)\circ'_px)_{\; I'}\circ' \lambda(S)^{r+1}\circ'_{J'} y\big)}_{(i)}-\underbrace{(-1)^{|x|} \Res \big( x_{\; I'}\circ' \lambda(S)^{r+1} \circ'_{J'} (\lambda(M)\circ'_p y)\big)}_{(ii)}\\
& +\underbrace{(-1)^{|x||y|} \Res \big( (\lambda(M) \circ'_q y)_{\; K'} \circ' \lambda(S)^{r+1} \circ'_{L'} x\big)}_{(iii)}\\
&+\underbrace{(-1)^{|x||y|+|y|}  \Res \big(y_{\; K'} \circ' \lambda(S)^{r+1} \circ'_{L'} (\lambda(M) \circ'_q x)\big)}_{(iv)}
\end{split}
\end{align*}
The summation is over indices $p$, $I'$, $J'$ such that $p\notin I'$
and $p\notin J'$, and similarly for $q$, $K'$ and $L'$. Let us compare
these terms with those in $\{\iota x, y\}_{r+1}$, which are given by
terms of the form
\[ \Res\big( (\lambda(M)\circ'_{p} x)_{\; I'}\circ' \lambda(S)^{r+1}\circ'_{J'} y \big)-(-1)^{|x||y|+|y|}\Res\big( y\circ'_{K'} \lambda(S)^{r+1}\circ'_{L'}(\lambda(M)\circ'_{p} x)\big) .\]
There are two cases, depending on whether the index $p$ is in $I'$ or not. 
\begin{itemize}
\item if $p\notin I'$ then these terms cancel precisely the terms $(i)$ and $(iv)$.
\item if $p\in I'$ the extra terms are equal to
\[ \sum_{\substack{p=i'_\alpha,\\1\leq \alpha\leq r+1}} \lambda(\mathbb{M}) \circ'_{(i'_\alpha,j'_\alpha)} \Res\big( x_{\; I'\backslash i'_\alpha}\circ' \lambda(S)^{r}\circ'_{J'\backslash j'_\alpha} y \big).\]
\end{itemize}
A similar argument also works for the terms in
$(-1)^{|x|}\{x, \iota y\}_{r+1}$, when a part of the terms cancel the
terms $(ii)$ and $(iii)$.  The extra term is given by
\[ (-1)^{|x||y|}  \sum_{\substack{q=k'_\alpha,\\1\leq \alpha\leq r+1}} \lambda(\mathbb{M}) \circ'_{(k'_\alpha,l'_\alpha)} \Res\big( y_{\; K'\backslash k'_\alpha}\circ' \lambda(S)^{r}\circ'_{L'\backslash l'_\alpha} x \big).\]
Adding the above two extra terms gives precisely $(r+1)\big(
\Delta\{x,y\}_r-\{\Delta x, y\}_r-(-1)^{|x|}\{x,\Delta y\}
\big)$. This proves the Leibniz rule. 

 The Jacobi identity is proved diagrammatically as in the proof of
Theorem~\ref{thm:dgla1}.  We have 
\begin{align*}
\{x,\{y,z\}\} & = \begin{tikzpicture}[box/.style={draw,rounded corners,text width=.5cm,align=center},scale=0.3,baseline={(current bounding box.center)}]
\node[box] at (-2,0) (b) {$x$};
\node[box] at (2,0) (c) {$y$};
\node[box] at (6,0) (d) {$z$};
\draw [thick,->] (b) to [out=-90, in=-90] node[midway,below]{} (1.8,-1);
\draw [thick,->] (2.2,-1) to [out=-90, in=-90] node[midway,below]{} (d);
\end{tikzpicture} - \begin{tikzpicture}[box/.style={draw,rounded corners,text width=.5cm,align=center},scale=0.3,baseline={(current bounding box.center)}]
\node[box] at (-2,0) (b) {$x$};
\node[box] at (2,0) (c) {$y$};
\node[box] at (6,0) (d) {$z$};
\draw [thick,<-] (b) to [out=-90, in=-90] node[midway,below]{} (1.8,-1);
\draw [thick,->] (2.2,-1) to [out=-90, in=-90] node[midway,below]{} (d);\end{tikzpicture}
\\
&+ \begin{tikzpicture}[box/.style={draw,rounded corners,text width=.5cm,align=center},scale=0.3,baseline={(current bounding box.center)}]
\node[box] at (-2,0) (b) {$x$};
\node[box] at (2,0) (c) {$y$};
\node[box] at (6,0) (d) {$z$};
\draw [thick,->] (b) to [out=-90, in=-90] node[midway,below]{} (6.2,-1);
\draw [thick,->] (2.2,-1) to [out=-90, in=-90] node[midway,below]{} (d);
\end{tikzpicture}-\begin{tikzpicture}[box/.style={draw,rounded corners,text width=.5cm,align=center},scale=0.3,baseline={(current bounding box.center)}]
\node[box] at (-2,0) (b) {$x$};
\node[box] at (2,0) (c) {$y$};
\node[box] at (6,0) (d) {$z$};
\draw [thick,<-] (b) to [out=-90, in=-90] node[midway,below]{} (6.2,-1);
\draw [thick,->] (2.2,-1) to [out=-90, in=-90] node[midway,below]{} (d);
\end{tikzpicture}\\
&-\begin{tikzpicture}[box/.style={draw,rounded corners,text width=.5cm,align=center},scale=0.3,baseline={(current bounding box.center)}]
\node[box] at (-2,0) (b) {$x$};
\node[box] at (2,0) (c) {$y$};
\node[box] at (6,0) (d) {$z$};
\draw [thick,->] (b) to [out=-90, in=-90] node[midway,below]{} (1.8,-1);
\draw [thick,<-] (2.2,-1) to [out=-90, in=-90] node[midway,below]{} (d);
\end{tikzpicture} +
\begin{tikzpicture}[box/.style={draw,rounded corners,text width=.5cm,align=center},scale=0.3,baseline={(current bounding box.center)}]
\node[box] at (-2,0) (b) {$x$};
\node[box] at (2,0) (c) {$y$};
\node[box] at (6,0) (d) {$z$};
\draw [thick,<-] (b) to [out=-90, in=-90] node[midway,below]{} (1.8,-1);
\draw [thick,<-] (2.2,-1) to [out=-90, in=-90] node[midway,below]{} (d);\end{tikzpicture}\\
&- \begin{tikzpicture}[box/.style={draw,rounded corners,text width=.5cm,align=center},scale=0.3,baseline={(current bounding box.center)}]
\node[box] at (-2,0) (b) {$x$};
\node[box] at (2,0) (c) {$y$};
\node[box] at (6,0) (d) {$z$};
\draw [thick,->] (b) to [out=-90, in=-90] node[midway,below]{} (6.2,-1);
\draw [thick,<-] (2.2,-1) to [out=-90, in=-90] node[midway,below]{} (d);
\end{tikzpicture}+
\begin{tikzpicture}[box/.style={draw,rounded corners,text width=.5cm,align=center},scale=0.3,baseline={(current bounding box.center)}]
\node[box] at (-2,0) (b) {$x$};
\node[box] at (2,0) (c) {$y$};
\node[box] at (6,0) (d) {$z$};
\draw [thick,<-] (b) to [out=-90, in=-90] node[midway,below]{} (6.2,-1);
\draw [thick,<-] (2.2,-1) to [out=-90, in=-90] node[midway,below]{} (d);
\end{tikzpicture}\\
&+\begin{tikzpicture}[box/.style={draw,rounded corners,text width=.5cm,align=center},scale=0.3,baseline={(current bounding box.center)}]
\node[box] at (-2,0) (b) {$x$};
\node[box] at (2,0) (c) {$y$};
\node[box] at (6,0) (d) {$z$};
\draw [thick,<-] (-1.6,-.8) to [out=-90, in=-90] node[midway,below]{} (1.8,-1);
\draw [thick,<-] (-2.4,-.8) to [out=-90, in=-90] node[midway,below]{} (6.2,-1);
\draw [thick,<-] (2.2,-1) to [out=-90, in=-90] node[midway,below]{} (d);
\end{tikzpicture}-\begin{tikzpicture}[box/.style={draw,rounded corners,text width=.5cm,align=center},scale=0.3,baseline={(current bounding box.center)}]
\node[box] at (-2,0) (b) {$x$};
\node[box] at (2,0) (c) {$y$};
\node[box] at (6,0) (d) {$z$};
\draw [thick,->] (-1.6,-.8) to [out=-90, in=-90] node[midway,below]{} (1.8,-1);
\draw [thick,->] (-2.4,-.8) to [out=-90, in=-90] node[midway,below]{} (6.4,-.8);
\draw [thick,<-] (2.2,-1) to [out=-90, in=-90] node[midway,below]{} (d);
\end{tikzpicture}\\
&-\begin{tikzpicture}[box/.style={draw,rounded corners,text width=.5cm,align=center},scale=0.3,baseline={(current bounding box.center)}]
\node[box] at (-2,0) (b) {$x$};
\node[box] at (2,0) (c) {$y$};
\node[box] at (6,0) (d) {$z$};
\draw [thick,<-] (-1.6,-.8) to [out=-90, in=-90] node[midway,below]{} (1.8,-1);
\draw [thick,<-] (-2.4,-.8) to [out=-90, in=-90] node[midway,below]{} (6.4,-.8);
\draw [thick,->] (2.2,-1) to [out=-90, in=-90] node[midway,below]{} (d);
\end{tikzpicture}+\begin{tikzpicture}[box/.style={draw,rounded corners,text width=.5cm,align=center},scale=0.3,baseline={(current bounding box.center)}]
\node[box] at (-2,0) (b) {$x$};
\node[box] at (2,0) (c) {$y$};
\node[box] at (6,0) (d) {$z$};
\draw [thick,->] (-1.6,-.8) to [out=-90, in=-90] node[midway,below]{} (1.8,-1);
\draw [thick,->] (-2.4,-.8) to [out=-90, in=-90] node[midway,below]{} (6.4,-.8);
\draw [thick,->] (2.2,-1) to [out=-90, in=-90] node[midway,below]{} (d);
\end{tikzpicture}.
\end{align*}
The signs arise from the number of left-to-right sewing
operations. Each arrow in the above diagram indicates a twisted sewing
operation of (possibly multiple) outputs with inputs. The last four
terms in the above sum are due to the fact that when performing the
second bracket operation, it is possible to choose inputs/outputs in
both $y$ and $z$.  Similarly, writing the terms in $\{\{x,y\},z\}$ and
$(-1)^{|x||y|}\{y,\{x,z\}\}$ all the $36$ terms cancel
out.
\end{Proof}

\paragraph{{\bf An example of the second construction of a DGLA.}} We
illustrate the second construction of a DGLA in the case of an
$\Annu$-algebra structure on a dg vector space $(V, b)$. We continue to
use the notations from~(\ref{para:example1}).

As a graded vector space we have
\[ \hg^V =\big( \bigoplus_{k\geq 1, l}
  \Hom(\wedge^k V_+, \Sym^l V_-)[2-k]\big)\series{\hbar}. \] 
The self-sewing operator $\Delta$ is only defined on
the output part of these $\Hom$'s, i.e. for
$\phi\in \Hom(\Sym( V_+[1]), \Sym^l V_-)$ we set
\[ (\Delta \phi) (x_1\cdots x_k)= \Delta\big(
  \phi(x_1\cdots x_k)\big),\] 
where on the right hand side the operator $\Delta$ acts by the formula
in~(\ref{para:example1}) applied to symmetric tensors.

We next consider the Lie bracket
\[ \{-,-\}_\hbar=\displaystyle\sum_{r\geq 1}
  \frac{1}{r!}\{-,-\}_r\cdot\hbar^{r-1}. \]
Denote by $\theta: V_-\ra V_+$ the map defined by 
\[ \theta( x_0+x_{-1}u^{-1}+\cdots)= uBx_0\in V_+. \] 
Since $\theta$ is odd (it has degree one), for each $r\geq 1$ we obtain
induced maps
\[ \theta^{(r)} : \Sym^r V_-  \ra \wedge^r V_+ \]
given by
\[ \theta^{(r)} (x_1\cdots x_r) = \theta(x_1)\wedge\cdots\wedge
  \theta(x_r). \]
Both the symmetric algebra $\Sym V_-$ and the exterior algebra
$\wedge V_+$ are cocommutative coalgebras with respect to the shuffle
coproduct.  The $r$-th Lie bracket of elements
$\phi \in \Hom(\wedge^{k} V_+, \Sym^l V_-)$ and
$\psi \in \Hom(\wedge^{k'} V_+, \Sym^{l'} V_-)$ lies inside
$\Hom(\wedge^{k+k'-r} V_+, \Sym^{l+l'-r} V_-)$.  It is given by
$\phi*\psi - (-1)^{|\phi||\psi|}\psi*\phi$ where $\phi*\psi$ is
defined as the composition
\begin{align*}
  \wedge^{k+k'-r}V_+ & \lra \wedge^{k'} V_+ \otimes \wedge^{k-r} V_+ \stackrel{\psi\otimes \id}{\longrightarrow} \Sym^{l'} V_- \otimes  \wedge^{k-r} V_+ \lra\\
  \Sym^r V_- \otimes & \Sym^{l'-r} V_- \otimes  \wedge^{k-r} V_+\stackrel{\theta^{(r)}\otimes\id}{\longrightarrow} \wedge^r V_+ \otimes \Sym^{l'-r} V_- \otimes  \wedge^{k-r} V_+ \lra\\
  \stackrel{\cong}{\longrightarrow} \wedge^r V_+ & \otimes\wedge^{k-r} V_+ \otimes \Sym^{l'-r} V_- \lra \wedge^k V_+ \otimes \Sym^{l'-r} V_- \lra\\
  \stackrel{\phi\otimes \id}{\longrightarrow}  \Sym^l & V_-   \otimes \Sym^{l'-r} V_- \lra \Sym^{l+l'-r} V_-.
\end{align*}

\paragraph{{\bf A map between the two constructions of DGLAs.}}
\label{para:iota}
Let $\cP$ be a PROP with $\Annu$ action given by
$\lambda: \Annu\ra \cP$.  By the above constructions we obtain two
DGLAs $\mathfrak{g}^\cP$ and $\widehat{\mathfrak{g}}^{\cP^+}$.  Denote
the part of the DGLA $\mathfrak{g}^\cP$ with positive number of
outputs by
\[ \big( \mathfrak{g}^{\cP} \big)^+= \bigoplus_{n\geq 1}
  \cP(0,n)_\hS[1]\series{\hbar}.\] 
Note that this is not a sub-DGLA of $\mathfrak{g}^\cP$,
but rather a quotient. Indeed, there exists a short exact sequence
 \[ 0\ra \cP(0,0)_\hS[1]\series{\hbar} \ra \mathfrak{g}^\cP \ra
   \big( \mathfrak{g}^{\cP} \big)^+ \ra 0\] 
which is a central extension of DGLAs. The relationship between the
constructions of Theorems~\ref{thm:dgla1} and~\ref{thm:dgla2} is
that there exists a natural morphism of DGLAs
\[ \iota: \big(\mathfrak{g}^\cP\big)^+ \ra \widehat{\mathfrak{g}}^{\cP^+}\]
which maps $x\in \cP(0,n)_\hS$ to 
\[ \iota(x) = \sum_{1\leq i\leq n} \lambda(M)\circ'_i x \in
  \cP(1,n-1)_\hS.\]
By construction, the image of this map is in the kernel of the
operator $\iota : P(1,n-1)_\hS \to P(2,n-2)_\hS$.  This is why the map
$\iota : (\g^\cP)^+ \to \hg^{\cP^+}$ is a cochain map.  It is a map of
Lie algebras for the simple reason that, on the subspace of elements
of $\hg^\cP$ with only one input, the only possible bracket involves
connecting only one input to one output.

For the examples of interest in this paper the map $\iota$ will often
be a quasi-isomorphism.

%% file: vertices.tex
\section{The Sen-Zwiebach DGLA and string vertices}
\label{sec:vertices}

We apply the constructions of the previous section to the dg PROP of
singular chains on the moduli spaces $C_*(M_{g,k,l}^\fr)$ of Riemann
surfaces with distinguished, framed inputs and outputs.  We find that
the map $\iota$ between the two types of DGLAs constructed before is a
quasi-isomorphism; since in the first one the Maurer-Cartan equation 
has solution that is unique up to homotopy, it follows that the same is true
for the second one.  The components of this unique solution are the
homological string vertices.  The same result holds for a combinatorial
version of the second type of DGLA above, constructed from ribbon
graphs. 

As stated in the introduction, $C_*(-)$ will denote the normalized
singular chain functor from topological spaces to dg vector spaces.

\paragraph{{\bf Moduli spaces of curves with framed markings.}} 
For $g\geq 0$ and $k,l\geq 0$ satisfying $2g-2+k+l>0$ we consider the
coarse moduli space $M_{g,k,l}^\fr$ of Riemann surfaces
of genus $g$ with $k+l$ marked and framed points, $k$ of which are
designated as inputs, and the rest as outputs.  Explicitly, a point in
$M_{g,k,l}^\fr$ is given by the data $(\Sigma, p_1,\ldots,p_k,q_1,\ldots,q_l,\phi_1,\ldots,\phi_k,\psi_1,\ldots,\psi_l)$ 
where
\begin{itemize}
\item $\Sigma$ is a smooth Riemann surface of genus $g$,
\item the $p$'s and the $q$'s are distinct marked points on $\Sigma$,
\item the $\phi$'s and the $\psi$'s are framings around the marked points, given by biholomorphic maps
\begin{align*}
\phi_i: &\, \mathbb{D}^2 \ra U_{\epsilon_i}(p_i),\;\; 1\leq i\leq k \\
\psi_j: &\, \mathbb{D}^2 \ra U_{\epsilon_j}(q_j), \;\; 1\leq j\leq l
\end{align*}
where $\mathbb{D}^2$ stands for the unit disk in $\mathbb{C}$, and
$U_\epsilon(x)$ is a radius $\epsilon$ disk of the marked point $x$ on
$\Sigma$, under, say, the unique hyperbolic metric on
$\Sigma\setminus \{p_1,\ldots,p_k,q_1,\ldots,q_l\}$.  We require that the
biholomorphic maps $\phi$'s and $\psi$'s extend to an open
neighborhood of $\mathbb{D}^2$.  We also require that the closures of all the framed disks be disjoint, i.e. for any two distinct points $x, y\in \{ p_1,\ldots,p_k,q_1,\ldots,q_l\}$ we have
\[\overline{U_{\epsilon_i}(x)} \cap  \overline{U_{\epsilon_j}(y)}=\emptyset.\]
\end{itemize}

We note that the fibration $M_{g,k,l}^\fr\ra M_{g,k,l}$ which
forgets the framing is a homotopy $(S^1)^{k+l}$-bundle.

\paragraph
\label{subsec:compo}
There is a sewing map 
\[ M_{h,l,m}^\fr \times M_{g,k,l}^\fr \ra M_{g+h+l-1,k,m}^\fr \]
defined by first removing the $l$ open disks from each framed
surface, and then sewing the boundary circles using the
identification $zw=1$. Here $z\in \partial \mathbb{D}^2$ and
$w\in \partial \mathbb{D}^2$ are coordinates on the two boundary
circles to be sewed together. Applying the functor $C_*(-)$ yields a
composition map 
\begin{equation*}
 \circ: C_*(M_{h,l,m}^\fr) \otimes C_*(M_{g,k,l}^\fr)
 \ra C_*(M_{g+h+l-1,k,m}^\fr). 
 \end{equation*} 

\paragraph{{\bf The PROP $\cS$.}} The collection of $\Sigma$-bimodules
$\{C_*(M_{g,k,l}^\fr)\}$ with composition defined above does not form
a PROP because we only allow connected surfaces in $M_{g,k,l}^\fr$. To
overcome this we set $\cR$ to be the free dg PROP generated by
$\{C_*(M_{g,k,l}^\fr)\}$ modulo the relations defined by the
composition maps $\circ$ in~(\ref{subsec:compo}). A geometric way to directly
construct $\cR$ would be to allow disconnected Riemann surfaces,
but we shall not need this description here.

We will be most interested in the PROP $\cS$ obtained by adjoining the
PROP $\Annu$ to $\cR$.  Conceptually this is equivalent to also
allowing the surfaces (annuli) in $M_{0,1,1}^\fr$ and $M_{0,2,0}^\fr$
-- in the framed setting these are stable (unlike in the unframed case).
To be precise we define $\cS$ to be the free PROP generated by
$\{C_*(M_{g,k,l}^\fr)\}$ and $S\in \Annu(1,1), M\in \Annu(2,0)$ modulo
the relations below.
\begin{itemize}
\item[--] The relations in the definition of $\Annu$, see~(\ref{para:annuli}).
\item[--] The relations in the definition of $\cR$, see~(\ref{subsec:compo}).
\item[--] For a $d$-chain $\sigma\in C_d(M_{g,k,l}^\fr)$ and an index
  $1\leq i\leq k$, define the composition $\sigma {\,}_{i} \circ S$
  to be (the prism decomposition of) the $(d+1)$-chain
  \[ S^1\times \Delta^d \stackrel{\id\times \sigma}{\longrightarrow}
    S^1\times M_{g,k,l}^\fr \ra M_{g,k,l}^\fr,\] 
  where the second arrow is the action of $S^1$ on the framing at the
  $i$-th marked point by precomposing a framing with the rotation map
  $e^{i\theta}:\mathbb{D}^2\ra \mathbb{D}^2$.  The
  composition on the other side $S\circ_j \sigma, \; (1\leq j\leq l)$
  is defined in the same way.
\item[--] For a $d$-chain $\sigma$ as above and a pair of indices
  $1\leq i<j\leq l$, define $M\circ_{(i,j)} \sigma$ as the
  composition
  \[ \Delta^d\stackrel{\sigma}{\longrightarrow} M^\fr_{g,k,l} \ra
    M_{g+1,k,l-2}^\fr.\] 
  Here the second map performs framed sewing at the two marked points $q_i$
  and $q_j$, i.e., we first remove the two disks $U_{\epsilon_i}(q_i)$,
  $U_{\epsilon_j}(q_j)$, and then we identify the two boundaries by $z=w$
  using the two framings $\psi_i$ and $\psi_j$.
\end{itemize}

\paragraph{{\bf Sen-Zwiebach's DGLA.}} The PROP $\Annu$ obviously sits
inside $\cS$, giving rise to an $\Annu$-action on $\cS$.
Theorem~\ref{thm:dgla1} constructs a DGLA 
\[ \g^\cS= \left ( \bigoplus_n \cS(0,n)_\hS[1]\right )\series{\hbar}.\] 
Inside $\cS(0,n)$ we have a subspace $C_*(M^\fr_{g,0,n})$ consisting
of connected Riemann surfaces.  We put these subspaces together to
form
\[ \g = \left (\bigoplus_{g,n} C_*(M^\fr_{g,0,n})_\hS[1]\right)\series{\hbar}.\]  
Since sewing operations with $S$ and $M$ preserves connectedness the
subspace $\g \subset \g^\cS$ is a sub-DGLA.

It is convenient to adjoin to $\g$ another formal variable $\lambda$
of degree $-2$.  We will use the notation $\g$ for the resulting DGLA
\[ \g= \left ( \bigoplus_{g,n} C_*(M^\fr_{g,0,n})_\hS [1]\right )
  \series{\hbar,\lambda} .\]

\begin{Definition}
\label{defi:vertices}
A degree $-1$ element $\cV\in \g$ of the form
\[ \cV=\sum_{g,n}\cV_{g,n}\hbar^g\lambda^{2g-2+n}, \;\; \cV_{g,n}\in
  C_*(M^\fr_{g,0,n})_\hS\] 
is called a homological string vertex if it satisfies the following properties:
\begin{enumerate}
\item $\cV$ is a Maurer-Cartan element of the DGLA $\g$, i.e., in
  terms of the components $\{\cV_{g,n}\}$ we have
\[ (\partial+uB) \cV_{g,n}+\Delta \cV_{g-1,n+2} +
  \frac{1}{2}\sum_{\substack{g_1+g_2=g\\n_1+n_2=n+2}}
  \{\cV_{g_1,n_1},\cV_{g_2,n_2}\} = 0. \]
\item $\cV_{0,3}=\frac{1}{6}\cdot {{\sf pt}}$ with ${{\sf pt}} \in
  C_0( M^\fr_{0,0,3})_\hS$ representing a point class. 
\end{enumerate}
\end{Definition}
\bigskip
\medskip

\noindent
A fundamental result about string vertices is the following theorem from~\cite{Cos}.
\vspace{-0.5mm}

\begin{Theorem}
\label{thm:costello}
The string vertex exists and is unique up to homotopy (i.e., gauge equivalence
between Maurer-Cartan elements). 
\end{Theorem}
\vspace{2mm}

\begin{Proof}
  The proof is a standard argument in deformation theory.  The
  obstruction to existence lies in $H_{6g-7+2n}(M^\fr_{g,0,n})_\hS$
  while the deformation space is given by
  $H_{6g-6+2n}(M^\fr_{g,0,n})_\hS$. Both homology groups are known to
  vanish. We refer to~\cite[Section 9]{Cos} for a detailed proof.
\end{Proof}

\paragraph{{\bf The Koszul resolution of $\g$.}} 
Let $\cS^+$ denote the sub-PROP of $\cS$ generated by chains on moduli
spaces of Riemann surfaces with positive number of inputs.  The
construction of Theorem~\ref{thm:dgla2}, applied to the morphism
$i: \Annu\ra \cS^+$, yields a second DGLA
\[\hg^{\cS^+}= \left (\bigoplus_{k\geq 1, l} \cS^+(k,l)_\hS[2-k]\right )\series{\hbar}.\]
Again, the subspace of connected stable surfaces forms a sub-DGLA to
which we adjoin the formal variable $\lambda$ to get 
\[ \hg= \left (\bigoplus_{g\geq 0, k\geq 1, l}
    C_*(M_{g,k,l}^\fr)_\hS[2-k]\right) \series{\hbar,\lambda}.\]
The canonical map $\iota: (\g^\cS)^+ \ra \hg^{\cS^+}$
of~(\ref{para:iota}) also restricts to give a morphism of DGLAs
$\iota: \g^+ \ra \hg$.

\begin{Proposition}
\label{prop:iota-qism}
The map $\iota: \g^+ \ra \hg$ is a quasi-isomorphism of DGLAs.
\end{Proposition}
\vspace{2mm}

\begin{Proof}
  Observe that for fixed $g\geq 0$ and $n>0$ the various moduli spaces
  $M^\fr_{g,k,l}$ with $k+l=n$ are all isomorphic. Then the result
  follows from the algebraic fact that if $W$ is any
  $\Sigma_n$-representation, then the associated Koszul complex
  \[ 0\ra W_{\Sigma_n} \ra W_{\Sigma_{1,n-1}} \ra
    W_{\Sigma_{2,n-2}}\ra \cdots\ra W_{\Sigma_{n,0}}\ra 0\] 
  is exact.  Here $W_{\Sigma_{k,l}}$ denotes the space of coinvariants
  of $W$ under the action of $\Sigma_k \times \Sigma_l$, where the
  first group acts via the sign representation $\sgn_k$ and the
  second one acts via the trivial representation.
\end{Proof}

\paragraph{{\bf The combinatorial version of $\hg$.}} 
There exists a version $\cS^{\comb,+}$ of $\cS^+$ constructed in terms
of ribbon graphs with framed inputs and outputs (also sometimes called
fat graphs with black and white vertices).  We will not give the
details of the construction of the PROP $\cS^{\comb, +}$ here since a
complete description of it will appear in the upcoming
work~\cite{CalChe}, following ideas of
Kontsevich-Soibelman~\cite{KonSoi} and Wahl-Westerland~\cite{WahWes}.
It suffices to say that the chain complexes $\cS^{\comb,+}(k,l)$ split
up as direct sums of complexes indexed by a genus $g\geq 0$,
\[ \cS^{\comb,+}(k,l) = \bigoplus_{g\geq 0} \cS^{\comb,+}_g(k,l). \]
Each summand $\cS^{\comb,+}_g(k,l)$ is quasi-isomorphic to
$C_*(M_{g,k,l}^\fr)$, and will be denoted by
$C_*^\comb(M_{g,k,l}^\fr)$.  A basis for $C_*^\comb(M_{g,k,l}^\fr)$
consists of isomorphism classes of framed ribbon graphs of genus $g$
with $k$ faces and $l$ white vertices.

The main results that we will need about $\cS^{\comb,+}$ are
summarized in the following theorems.

\begin{Theorem}
  There is a sub-PROP of $\cS^{\comb,+}$ equivalent to $\Annu$.
\end{Theorem}
\medskip

\begin{Proof}
  The operations $S\in \Annu(1,1)$ and $M\in \Annu(2,0)$ correspond to
  the following ribbon graphs
  \[ S= \begin{tikzpicture}[baseline={([yshift=-1.2ex]current bounding box.center)},scale=0.3]
\draw (0,0) node[cross=2pt,label=above:{}] {};
\draw[thick]  (0,0) to (1.2,0);
\draw (1.4,0) circle (.2);
\draw [ultra thick] (1.6,0) to (2.4,0);
\end{tikzpicture}, \;\;\;\;\;\;\;\;\;
M= \begin{tikzpicture}[baseline={([yshift=-0.4ex]current bounding box.center)},scale=0.3]
\draw [thick] (0,2) circle [radius=2];
\draw [thick] (-2,2) to (-0.6,2);
\draw [thick] (-0.8,2.2) to (-0.4,1.8);
\draw [thick] (-0.8, 1.8) to (-0.4, 2.2);
\draw [thick] (2,2) to (3.4,2);
\draw [thick] (3.2,2.2) to (3.6,1.8);
\draw [thick] (3.2, 1.8) to (3.6, 2.2);
\end{tikzpicture} .\]
The thick leaf in the graph $S$ is a marked leaf of the white vertex (see~\cite{WahWes}) which indicates a generic framing for the marked point corresponding to the white vertex.
\end{Proof}

\begin{Theorem}
  The PROP $\cS^{\comb,+}$ is quasi-equivalent to $\cS$.
\end{Theorem}
\bigskip

\begin{Proof}
  This follows from work of Egas~\cite{Ega}.
\end{Proof}

\begin{Theorem}
  Let $A$ be a cyclic $A_\infty$-algebra of Calabi-Yau dimension $d$.
  Then $A$ gives rise to a 2d TFT structure on the dg-vector space
  $L=CC_*(A)[d]$ of shifted Hochschild chains of $A$, i.e., a PROP map
  \[ \rho^A: \cS^{\comb,+} \lra \PEnd(L). \]
\end{Theorem}
\vspace*{-2mm}

\begin{Proof}
  This constructions was sketched in~\cite[Section 11.6]{KonSoi}.  Complete details of this result will appear
  in~\cite{CalChe}.
\end{Proof}

\begin{Theorem}
  The ribbon graphs that form a basis of $C_*^\comb(M_{g,k,l}^\fr)$
  and the operators $b$, $B$, $\iota$, $\Delta$, $\{-,-\}_i$
  can be algorithmically computed in explicit form.
\end{Theorem}
\bigskip

\begin{Proof}
  This follows from the explicit description of $\cS^{\comb,+}$
  in~\cite{CalChe}.
\end{Proof}

\paragraph
To summarize, we have the following diagram:
\[\begin{tikzpicture}
\node at (0,1.8) {$\Annu$};
\draw [thick,->] (.3,1.5) to (2,0);
\draw [thick,->] (0,1.5) to (0,0);
\draw [thick,->] (-0.3,1.5) to (-2,0);
\node at (-2, -.2) {$\cS$};
\node at (0,-.2) {$\cS^+$};
\node at (2.4,-.2) {$\cS^{+,{\comb}}$};
\draw [thick,->] (-.3,-.2) to (-1.8,-.2);
\draw [thick,<-]  (.3,-.2) to (1.7,-.2);
\node at (1, 0) {$\cong$};
\end{tikzpicture}\]
This yields a roof diagram of quasi-isomorphisms of DGLAs:
\[\begin{tikzpicture}
\node at (-3,-1) {$\displaystyle{\hg=\left (\bigoplus_{g, k\geq 1, l}
    C_*(M_{g,k,l}^\fr)_\hS[2-k]\right )\series{\hbar}}$};
\draw [thick,->] (0,1.5) to (-2,0);
\node at (1, .7) {$\cong$ (by Proposition~\ref{prop:iota-qism})};
\draw [thick,<-] (-2,-1.5) to (0, -3);
\node at (1.3, -2.2) {$\cong$ (since $\cS^{+,\comb}\cong \cS^+$)};
\node at (0,1.8)
{$\g=\displaystyle{\left
      (\bigoplus_{g,n>0}C_*(M^\fr_{g,0,n})_\hS[1]\right )
    \series{\hbar}}$}; 
\node at (0,-4) {$\displaystyle{\hg^{\comb}=\left (\bigoplus_{g, k\geq
        1, l} C^{\comb}_*(M_{g,k,l}^\fr)_\hS[2-k] \right ) \series{\hbar}}$};
\end{tikzpicture}\]
By the homotopy invariance of Maurer-Cartan moduli space,
Theorem~\ref{thm:costello} implies the following result. 

\begin{Theorem}~\label{thm:comb-vertex}
In the DGLA
 \[ \hg^{{\comb}}= \left ( \bigoplus_{g, k\geq 1, l}
   C^{\comb}_*(M_{g,k,l}^\fr)_\hS[2-k] \right )
 \series{\hbar,\lambda} \]
there exists a degree $-1$ element $\hcV^\comb$, unique up to homotopy, of the
form
\[\hcV^\comb=\sum_{g, k\geq 1, l} \hcV^\comb_{g,k,l}\,\hbar^g
\lambda^{2g-2+k+l}, \]
such that the following conditions are satisfied:
\begin{enumerate}
\item $\hcV^\comb$ is a Maurer-Cartan element of $\hg^{{\comb}}$,
  i.e., $\hcV^\comb$ satisfies the following equations for each triple
  $(g,k\geq 1,l)$:
\[ (\partial'+uB')
  \hcV^\comb_{g,k,l}+\iota\hcV^\comb_{g,k-1,l+1}+\Delta\hcV^\comb_{g-1,k,l+2}+\frac{1}{2}\sum\frac{1}{r!}\{\hcV^\comb_{g_1,k_1,l_1},\hcV^\comb_{g_2,k_2,l_2}\}_r=0.\]
The last sum is over all $r\geq 1$ and all $(g_1, g_2, k_1, k_2, l_1,
l_2)$ such that
\begin{align*}
  g_1+g_2 +r -1 & = g\\
  k_1+k_2 -r & = k\\
  l_1+l_2 -r & = l
\end{align*}

\item $\hcV^\comb_{0,1,2}=\frac{1}{2}\cdot {{\sf pt}}$ with ${{\sf
      pt}}$ representing a point class in $C^{\comb}_0(
  M^\fr_{0,1,2})_\hS$, i.e., the ribbon graph expression
 \begin{align*}
\hcV^\comb_{0,1,2}& =\frac{1}{2}\;\;\;
  \begin{tikzpicture}[baseline={([yshift=-2ex]current bounding
      box.center)},scale=0.3] 
\draw [thick] (0,0) to (0,2);
\draw [thick] (-0.2, 1.8) to (0.2, 2.2);
\draw [thick] (0.2, 1.8) to (-0.2, 2.2);
\draw [thick] (0,0) to (-2,0);
\draw [thick] (0,0) to (2,0);
\draw [thick] (-2.2,0) circle [radius=0.2];
\draw [thick] (2.2,0) circle [radius=0.2];
\end{tikzpicture}
\end{align*}
\end{enumerate}
\end{Theorem}
\medskip

\begin{Definition}
  We shall refer to the $\bbQ$-linear combinations of ribbon graphs
  $\hcV^\comb_{g,k,l}$ as {{\sl combinatorial string vertices}}.
\end{Definition}
\medskip

\paragraph{{\bf Explicit formulas for some combinatorial string
    vertices.}} 
We present explicit formulas for the first few combinatorial string
vertices, ignoring orientations of ribbon graphs and associated
signs.  Our conventions are that a factor of $u^{-k}$ associated to a
boundary component of a diagram is related to the class $\psi^{k-1}$.
If we do not indicate a power of $u$, then we mean $u^{-1}$.
\begin{align*}
\hcV^\comb_{0,1,2}& =\frac{1}{2}\;\;\;
  \begin{tikzpicture}[baseline={([yshift=-1.2ex]current bounding
      box.center)},scale=0.3] 
\draw [thick] (0,0) to (0,2);
\draw [thick] (-0.2, 1.8) to (0.2, 2.2);
\draw [thick] (0.2, 1.8) to (-0.2, 2.2);
\draw [thick] (0,0) to (-2,0);
\draw [thick] (0,0) to (2,0);
\draw [thick] (-2.2,0) circle [radius=0.2];
\draw [thick] (2.2,0) circle [radius=0.2];
\end{tikzpicture}
\end{align*}

\[ \hcV^\comb_{1,1,0}=\frac{1}{24}\;\;\; \begin{tikzpicture}[baseline={([yshift=-.4ex]current bounding box.center)},scale=0.3]
\draw [thick] (0,2) circle [radius=2];
\draw [thick] (-2,2) to (-0.6,2);
\draw [thick] (-0.8,2.2) to (-0.4,1.8);
\draw [thick] (-0.8, 1.8) to (-0.4, 2.2);
\draw [thick] (-1.4142, 0.5858) to [out=40, in=140] (1, 0.5);
\draw [thick] (1.25, 0.3) to [out=-45, in=225] (2, 0.2);
\draw [thick] (2,0.2) to [out=45, in=-50] (1.732, 1);
\node at (0.7, 2.2) {$\scriptstyle{u^{-2}}$};
\end{tikzpicture}
\;\;+\;\;\frac{1}{4}\;\;\; 
\begin{tikzpicture}[baseline={([yshift=-.4ex]current bounding box.center)},scale=0.3]
\draw [thick] (0,2) circle [radius=2];
\draw [thick] (0,0) to (0,1.4);
\draw [thick] (-0.2, 1.2) to (0.2, 1.6);
\draw [thick] (-0.2, 1.6) to (0.2, 1.2);
\draw [thick] (0,0) to [out=80, in=180] (0.5, 1);
\draw [thick] (0.5,1) to [out=0, in=100] (0.9, 0.4);
\draw [thick] (0,0) to [out=-80, in=180] (0.5, -1);
\draw [thick] (0.5, -1) to [out=0, in=-100] (0.9, 0);
\end{tikzpicture}\] 

\[ \hcV^\comb_{0,2,1}= \frac{1}{2} \begin{tikzpicture}[baseline={([yshift=-.4ex]current bounding box.center)},scale=0.3]
\draw (0,0) node[cross=2pt,label=above:{}] {};
\draw[ultra thick]  (0,0) to (0.8,0);
\draw (1,0) circle (.2);
\draw (5.2,0) node[cross=2pt,label=above:{}] {};
\draw [thick] (5.2,0) to (6.2,0);
\draw [thick] (1.2,0) to (2.2,0);
\draw [thick] (4.2,0) circle [radius=2];
\end{tikzpicture}+\frac{1}{2} \begin{tikzpicture}[baseline={([yshift=-.4ex]current bounding box.center)},scale=0.3]
\draw (2,0) node[cross=2pt,label=above:{}] {};
\draw (6.2,0) node[cross=2pt,label=above:{}] {};
\draw [thick] (6.2,0) to (7.2,0);
\draw [thick] (2.1,0) to (3.2,0);
\draw [thick] (5.2,0) + (-85:2) arc(-85:265:2);
\draw [thick] (5.2,-2) circle (.2);
\draw [ultra thick,domain=230:265] plot ({5.2+2*cos(\x)}, {2*sin(\x)});
\end{tikzpicture}+\begin{tikzpicture}[baseline={([yshift=-.4ex]current bounding box.center)},scale=0.3]
\draw (2,0) node[cross=2pt,label=above:{}] {};
\draw (6.2,0) node[cross=2pt,label=above:{}] {};
\draw [thick] (6.2,0) to (7.2,0);
\draw [thick] (2.1,0) to (3.2,0);
\draw [thick] (5.2,0) circle [radius=2];
\draw (1.8,-1.6) circle (.2);
\draw [thick] (3.2,0) to (1.9,-1.4);
\end{tikzpicture}\]

\[\hcV^\comb_{0,3,0}= \frac{1}{2} \begin{tikzpicture}[baseline={([yshift=-.4ex]current bounding box.center)},scale=0.3]
\draw (2,0) node[cross=2pt,label=above:{}] {};
\draw (5.2,1) node[cross=2pt,label=above:{}] {};
\draw (5.2,-1) node[cross=2pt,label=above:{}] {};
\draw [thick] (5.2,-1) to (7.2,0);
\draw [thick] (5.2,1) to (5.2,2);
\draw [thick] (3.2,0) to (7.2,0);
\draw [thick] (2.1,0) to (3.2,0);
\draw [thick] (5.2,0) circle [radius=2];
\end{tikzpicture}+
\frac{1}{2} \begin{tikzpicture}[baseline={([yshift=-.4ex]current bounding box.center)},scale=0.3]
\draw [thick] (0,0) circle [radius=1.5];
\draw [thick] (3,0) circle [radius=1.5];
\draw (-.5,0) node[cross=2pt,label=above:{}] {};
\draw (2.5,0) node[cross=2pt,label=above:{}] {};
\draw (5.5,0) node[cross=2pt,label=above:{}] {};
\draw [thick] (5.5,0) to (4.5,0);
\draw [thick] (2.5,0) to (1.5,0);
\draw [thick] (-.5,0) to (-1.5,0);
\end{tikzpicture}+
\frac{1}{2} \begin{tikzpicture}[baseline={([yshift=-1ex]current bounding box.center)},scale=0.3]
\draw [thick] (0,0) circle [radius=1.5];
\draw [thick] (3,0) circle [radius=1.5];
\draw (-.5,0) node[cross=2pt,label=above:{}] {};
\draw (1.5,1.5) node[cross=2pt,label=above:{}] {};
\draw (3.5,0) node[cross=2pt,label=above:{}] {};
\draw [thick] (3.5,0) to (4.5,0);
\draw [thick] (1.5,1.5) to (1.5,0);
\draw [thick] (-.5,0) to (-1.5,0);
\end{tikzpicture}\]

\[\hcV^\comb_{0,1,3} =\frac{1}{2}\;
  \begin{tikzpicture}[baseline={([yshift=-1.2ex]current bounding
      box.center)},scale=0.3] 
\draw (0,1.5) node[cross=2pt,label=above:{}] {};
\draw [thick] (0,0.2) to (0,1.5);
\draw [thick] (-0.2,0) to (-2,0);
\draw [thick] (0.2,0) to (2,0);
\draw [thick] (-2.2,0) circle [radius=0.2];
\draw [thick] (2.2,0) circle [radius=0.2];
\draw [thick] (0,0) circle [radius=0.2];
\draw [ultra thick] (.2,0) to (1,0);
\end{tikzpicture}+\frac{1}{2}\;
  \begin{tikzpicture}[baseline={([yshift=-1.2ex]current bounding
      box.center)},scale=0.3] 
\draw (1,1) node[cross=2pt,label=above:{}] {};
\draw [thick] (1,0) to (1,1);
\draw [thick] (-0.2,0) to (-2,0);
\draw [thick] (0.2,0) to (2,0);
\draw [thick] (-2.2,0) circle [radius=0.2];
\draw [thick] (2.2,0) circle [radius=0.2];
\draw [thick] (0,0) circle [radius=0.2];
\draw [ultra thick] (0,-.2) to (0,-1);\end{tikzpicture}
+\frac{1}{2}\;
  \begin{tikzpicture}[baseline={([yshift=-.4ex]current bounding
      box.center)},scale=0.3] 
\draw (1,1) node[cross=2pt,label=above:{}] {};
\draw [thick] (1,0) to (1,1);
\draw [thick] (0,0) to (-2,0);
\draw [thick] (0,0) to (2,0);
\draw [thick] (-2.2,0) circle [radius=0.2];
\draw [thick] (2.2,0) circle [radius=0.2];
\draw [thick] (0,-1) circle [radius=0.2];
\draw [thick] (0,0) to (0,-.8);
\node at (1.2,-1) {$\scriptstyle{u^{-2}}$};
\end{tikzpicture}+\frac{1}{6}\;
  \begin{tikzpicture}[baseline={([yshift=-2ex]current bounding
      box.center)},scale=0.3] 
    \draw (-1,1) node[cross=2pt,label=above:{}] {}; \draw [thick]
    (-1,0) to (-1,1); \draw [thick] (0,0) to (-2,0); \draw [thick]
    (0,0) to (2,0); \draw [thick] (-2.2,0) circle [radius=0.2]; \draw
    [thick] (2.2,0) circle [radius=0.2]; \draw [thick] (1,1) circle
    [radius=0.2]; \draw [thick] (1,0) to (1,.8); \node at (-1,1.5)
    {$\scriptstyle{u^{-2}}$};\end{tikzpicture}.\]
Note, for example, that the coefficients $1/24$ and $1/6$ in
$\hcV_{1,1,0}^\comb$ and $\hcV_{0,1,4}^\comb$, respectively,
correspond to the integrals of a $\psi$-class on $\bM_{1,1}$ and
$\bM_{0,4}$, respectively.  This will be explained in~\cite{CalTu}.

%% file: tcft.tex
\section{String functionals from topological field theories}
\label{sec:tcft}

In this section we consider the action of the combinatorial string
vertices constructed in the previous section on the Hochschild chains
of a cyclic $\cA_\infty$-algebra.  The resulting functionals
$\widehat{\beta}_{g,k,l}$ are natural homological analogues of the string
functionals introduced by Sen-Zwiebach in~\cite{SenZwi}.

In this section we assume that all $A_\infty$-algebras are defined
over a field $\bbk$ of characteristic zero.

\paragraph{{\bf Two dimensional TFTs of dimension $d$.}}  For a non-negative
integer $d\geq 0$ define a shifted version $\cS_{[d]}$ of the PROP
$\cS$ by performing the following degree shifts on the generators:
\begin{align*}
\begin{split} C_*(M_{g,k,l}^{{\sf fr}}) & \mapsto C_*(M_{g,k,l}^{{\sf
fr}})[d(2-2g-2k)],\\ S & \mapsto S,\\ M & \mapsto M[-2d].
\end{split}
\end{align*} In other words, after the shift $M$ has degree $-2d$,
while $S$ still has degree $1$.  The composition map $\circ$
from~(\ref{subsec:compo}) still has degree zero.  We have similar
shifted versions of $\cS^+$, $\cS^{+,\comb}$ which will be denoted by
$\cS^+_{[d]}$ and$\cS^{+,\comb}_{[d]}$, respectively.
\newpage

Following~\cite{Cos} we define:
\begin{itemize}
\item[--] A two dimensional TFT of dimension $d$ is a morphism of dg PROPs
\[ F: \cS_{[d]} \ra \PEnd(V).\]
\item[--] A positive boundary two dimensional TFT is a morphism of dg PROPs
\[ F: \cS^+_{[d]} \ra \PEnd(V).\]
\item[--] A combinatorial positive boundary two dimensional TFT is a morphism of dg
PROPs 
\[ F: \cS^{+,\comb}_{[d]} \ra \PEnd(V).\]
\end{itemize}

The following result of Kontsevich-Soibelman~\cite{KonSoi}
and~\cite{Cos1} provides us with non-trivial examples of combinatorial
positive boundary two dimensional TFTs.  Explicit formulas for the PROP action will
be given in~\cite{CalChe} following the ideas in~\cite{KonSoi}
and~\cite{WahWes}.

\begin{Theorem}
  \label{thm:tcft} Let $A$ be a cyclic unital $A_\infty$-algebra of
  Calabi-Yau dimension $d$. Then there exists a combinatorial positive
  boundary two dimensional TFT of dimension $d$ given by a dg PROP morphism
\[ \rho^A: \cS^{+,\comb}_{[d]} \ra \PEnd(L),\] whose underlying dg
vector space is the shifted reduced Hochschild chain complex
$L=CC_*(A)[d]$ of the $A_\infty$-algebra $A$.
\end{Theorem}
\vspace*{2mm}

\paragraph{{\bf The $\Annu$-algebra structure.}} An immediate
consequence of Theorem~\ref{thm:tcft} is the existence of an
$\Annu$-algebra structure on the dg vector space $L$. The degree one
operator $B =\rho^A(S): CC_*(A)[d] \ra CC_*(A)[d]$ is the Connes
cyclic differential. The degree $-2d$ operator
\[ \langle-,-\rangle_{{\sf Muk}} = \rho^A(M): CC_*(A)[d]\otimes
  CC_*(A)[d] \ra \mathbb{K}\]
is known as the Mukai pairing. It is a symmetric pairing of degree
$-2d$. Applying the constructions of Theorems~\ref{thm:dgla1}
and~\ref{thm:dgla2} we obtain two DGLAs naturally associated with the
$\Annu$-algebra $L$. Denote them by
\begin{align*} \h & = \left ((\Sym L_-)[1]\right )\series{\hbar,\lambda}\\
\hh & = \left( \bigoplus_{k\geq 1, l} \Hom^\cont (\wedge^k L_+, \Sym^l
      L_-)[2-k]\right)\series{\hbar,\lambda} 
\end{align*}
Here we continue to use the notations of~(\ref{para:example1}):
$L_+=uL\series{u}$ and $L_-=L[u^{-1}]$.  We also denote by
$\Hom^\cont$ the space of continuous homomorphisms in the $u$-adic
topology ($\hh$ is a Lie subalgebra of the one defined in
Theorem~\ref{thm:dgla2}). 

\begin{Lemma}
  \label{lem:gamma}
  Let $A$ be a smooth, cyclic, and unital $A_\infty$-algebra which is
  assumed to satisfy the Hodge-de Rham degeneration
  property~\footnote{If $A$ is $\Z$-graded the degeneration property
    is a consequence of $A$ being smooth by a result of
    Kaledin~\cite{Kal}.}. Then the morphism defined
  in~(\ref{para:iota})
  \[ \iota: \h^+= \left ( (\Sym^{\geq 1}
  L_-)[1]\right ) \series{\hbar,\lambda}\longrightarrow \hh \]
  is a quasi-isomorphism of DGLAs.
\end{Lemma}
\vspace{4mm}

\begin{Proof}
It suffices to prove that for each $n\geq 1$ the following sequence of
dg vector spaces is exact:
\[ 0\ra \Sym^n L_- \ra \Hom^\cont( L_+, \Sym^{n-1} L_-)\ra \cdots \ra
  \Hom^\cont( \wedge^n L_+, \mathbb{C}) \ra 0. \]
We argue that the total complex of the sequence has zero homology, by
considering the spectral sequence associated with the
$u$-filtration. Since $A$ satisfies the degeneration property the
$E^1$ page is given by
\[ 0\ra \Sym^n H_- \ra \Hom^\cont( H_+, \Sym^{n-1} H_-)\ra \cdots \ra
  \Hom^\cont( \wedge^n H_+, \mathbb{C}) \ra 0 \]
where $H=HH_*(A)[d]$ is the shifted Hochschild homology of $A$. Because
$A$ is assumed to be smooth, $H$ is finite dimensional. From this
finiteness condition we can deduce the exactness of the $E^1$-page as
follows. Observe that there is an isomorphism for each $k+l=n$
\[ I_{k,l}: \wedge^k H_- \otimes \Sym^l H_- \ra \Hom^\cont( \wedge^k
  H_+, \Sym^{l} H_-)\]
defined using the Mukai pairing on $H$. In the case of $k=1$,
$l=0$ this map is explicitly given by
\[ I_{1,0} (x) (y) = \res_{u=0} \langle x, y\rangle.\]
For general $k$ and $l$ it is defined similarly. The fact that
$I_{k,l}$ is an isomorphism follows from a result of
Shklyarov~\cite{Shk} proving that the categorical Mukai pairing is
non-degenerate when $A$ is smooth (and finite dimensional, which
follows from the cyclic property of $A$).  Putting the isomorphisms
$I_{k,l}$ together we conclude that the $E^1$ page is isomorphic to
\[ 0\ra \Sym^n H_- \ra H_-\otimes \Sym^{n-1} H_- \ra \wedge^2
H_-\otimes \Sym^{n-2} H_- \ra \cdots \wedge^n H_- \ra 0,\] which is
endowed with the usual Koszul differential.  The exactness
follows.
\end{Proof}

\paragraph{{\bf Pushing forward the string vertex.}}
\label{para:mc}
By Theorem~\ref{thm:tcft} there is a combinatorial positive boundary
two dimensional TFT given by a morphism of dg PROPs
\[ \rho^A: \cS^{+,\comb}_{[d]} \ra \PEnd\big(L\big).\] The PROP
$\Annu_{[d]}$ sits inside $\cS^{+,\comb}_{[d]}$ which gives us a
commutative diagram of PROPs
\[\begin{tikzpicture} \node at (0,1.8) {$\Annu_{[d]}$}; \draw
[thick,->] (.3,1.5) to (1.8,0); \draw [thick,->] (-0.3,1.5) to
(-1.8,.2); \node at (-2, -.2) {$\cS^{+,\comb}_{[d]}$}; \node at
(2.4,-.2) {$\PEnd\big(L\big).$}; \draw [thick,->] (-1.5,-.2) to
(1.7,-.2);
\end{tikzpicture}\]
Since the degree shifts in the definition of the shifted PROPS above
are all even, the action $\rho^A$ induces a morphism (which we still
denote by $\rho^A$) of $\Z/2\Z$-graded DGLAs
\[ \rho^A: \hg \ra \hh. \]
The combinatorial string vertex constructed in Theorem~\ref{thm:comb-vertex},
\[\hcV^\comb=\sum_{g,k\geq 1, l}
\hcV^\comb_{g,k,l}\,\hbar^g \lambda^{2g-2+k+l}\]
yields a Maurer-Cartan element of $ \hh$ of the form
\[\widehat{\alpha}^A= \rho^A(\hcV^\comb)=
\sum_{g, k\geq 1, l} \widehat{\alpha}^A_{g,k,l}\,\hbar^g
\lambda^{2g-2+k+l}\]
with
\[ \widehat{\alpha}^A_{g,k,l} = \rho^A(\hcV_{g,k,l}^\comb) \in
  \Hom^\cont(\wedge^{k} L_+, \Sym^l  L_-). \]
Note that by the uniqueness of $\hcV^\comb$ the
Maurer-Cartan element $\widehat{\alpha}^A\in \hh$ is
also uniquely defined up to gauge equivalence (though the
Maurer-Cartan equation in $\hh$ may have many non-equivalent
solutions!).

Applying to the tensors $\widehat{\alpha}^A_{g,k,l}$ the sign
corrections in~(\ref{subsec:beta}) we obtain tensors
\[ \widehat{\beta}^A_{g,k,l} \in \Hom^\cont\left(\Sym^k(L_+[1]),
    \Sym^l L_-\right ) \]
such that
\[ \widehat{\beta} = \sum_{g,k\geq 1, l}
  \widehat{\beta}^A_{g,k,l}\hbar^g\lambda^{2g-2+k+l} \]
satisfies the Maurer-Cartan equation in an appropriately modified
version of $\hh$.

Furthermore, by Lemma~\ref{lem:gamma} there exists
a unique (up to gauge equivalence) Maurer-Cartan element
$\beta^A\in \h^+$ such that
\begin{itemize}
\item $\beta^A= \sum_{g, n\geq 1}\beta^A_{g,n} \,\hbar^g\lambda^{2g+n-2}$ with
$\beta^A_{g,n}\in \Sym^n L_-$.
\item $\iota(\beta^A)$ is gauge equivalent to
$\widehat{\beta}^A=\rho^A\big( \hcV^\comb\big)$.
\end{itemize}
The elements $\beta^A$ and $\widehat{\beta}^A$ are homological
analogues of the string functionals defined by
Sen-Zwiebach~\cite{SenZwi}.  They will be used in~\cite{CalTu} to
define the categorical enumerative invariants.